\documentclass[a4paper]{amsart}

\title[Jump-Diffusions in Hilbert Spaces: Existence, Stability and Numerics]{Jump-Diffusions in Hilbert Spaces: Existence, Stability and Numerics}
\author{Damir Filipovi\'c \and Stefan Tappe \and Josef Teichmann}
\address{Vienna Institute of Finance, University of Vienna, and Vienna University of Economics and Business Administration, Heiligenst\"adter Strasse 46-48, A-1190 Wien, Austria; Vienna University of Technology, Department of Mathematical Methods in Economics, Wiedner Hauptstrasse 8--10, A-1040 Wien, Austria}
\email{stefan.tappe@vif.ac.at, damir.filipovic@vif.ac.at,\newline jteichma@fam.tuwien.ac.at}
\thanks{The first and second author gratefully acknowledge the support from WWTF (Vienna Science and Technology Fund). The third author gratefully acknowledges the support from the FWF-grant Y 328 (START prize from the Austrian Science Fund).}
\thanks{The authors thank Christa Cuchiero for discussions about $L^p$-estimates and Martin Hairer for discussions about the Sz\H{o}kefalvi-Nagy theorem. The authors thank an anonymous referee for very helpful comments and suggestions.}

\usepackage{amscd}
\usepackage{amsmath}
\usepackage{amssymb}
\usepackage{amsthm}
\usepackage{bbm}
\usepackage{stmaryrd}

\newif\ifpdf
\ifx\pdfoutput\undefined
   \pdffalse        
\else
   \pdfoutput=1     
   \pdftrue
\fi

\ifpdf
   \usepackage[pdftex]{graphicx}
\else
   \usepackage{graphicx}
\fi

\numberwithin{equation}{section} \swapnumbers

\newtheorem{satz}{Satz}[section]

\newtheorem{theorem}[satz]{Theorem}
\newtheorem{proposition}[satz]{Proposition}
\newtheorem{corollary}[satz]{Corollary}
\newtheorem{lemma}[satz]{Lemma}
\newtheorem{assumption}[satz]{Assumption}

\newtheorem{definition}[satz]{Definition}

\newtheorem{remark}[satz]{Remark}
\newtheorem{remarks}[satz]{Remarks}

\newtheorem{example}[satz]{Example}

\newcommand{\abs}[1]{\lvert #1 \rvert}

\newcommand{\Lip}{\operatorname{Lip}}
\newcommand{\R}{\mathbb{R}}
\newcommand{\N}{\mathbb{N}}

\newcommand{\f}[2]{\frac{#1}{#2}}

\begin{document}

\begin{abstract}
By means of an original approach, called ``method of the moving frame'', we establish existence, uniqueness and stability results for mild and weak solutions of stochastic partial differential equations (SPDEs) with path dependent
coefficients driven by an infinite dimensional Wiener process and a compensated Poisson random measure. Our approach is based on a time-dependent coordinate transform, which reduces a wide class of SPDEs to a class of simpler SDE problems. We try to present the most general results, which we can obtain in our setting, within a self-contained framework to demonstrate our approach in all details. Also several numerical approaches to SPDEs in the spirit of this setting are presented.
\end{abstract}

\maketitle\thispagestyle{empty}


\bigskip


\bigskip

\section{Introduction}

Stochastic partial differential equations (SPDEs) are usually considered as stochastic perturbations of partial differential equations (PDEs). More precisely, let $ H $ be a Hilbert space and $ A $ the generator of a strongly continuous semigroup $ S $ on $ H $, then
$$
\frac{d r_t}{dt} = A r_t + \alpha(r_t), \quad r_0 \in H
$$
describes a (semi-linear) PDE on the Hilbert space of states $H$ with linear generator $A$ and (non-linear) term $ \alpha: H \to H $. Solutions are usually defined in the mild or weak sense. A stochastic perturbation of this (semi-linear) PDE is given through a driving noise and (volatility) vector fields, for instance one can choose a one-dimensional Brownian motion $ W $ and $ \sigma : H \to H $ and consider
$$
d r_t = (A r_t + \alpha(r_t))dt + \sigma(r_t) dW_t, \quad r_0 \in H.
$$
Solution concepts, properties of solutions, manifold applications have been worked out in the most general cases, e.g., \cite{Da_Prato} in the case of Brownian noise or \cite{P-Z-book} in the case of L\'evy noises.

We suggest in this article a new approach to SPDEs, which works for
most of the SPDEs considered in the literature (namely those where
the semigroup is pseudo-contractive). The advantages are three-fold:
first one can consider most general noises with path-dependent
coefficients and derive existence, uniqueness and stability results
in an easy manner. Second the new approach easily leads to
(numerical) approximation schemes for SPDEs, third the approach
allows for rough path formulations (see \cite{Tei:08}) and therefore
for large deviation results, Freidlin-Wentzell type results, etc. In
this article we shall mainly address existence, uniqueness and
stability results for SPDEs with driving Poisson random measures and
general path-dependent coefficients. An outline of the basic
relation, namely short-time asymptotics, for high-order, weak or
strong numerical schemes is presented, too.

In our point of view SPDEs are considered as time-dependent
transformations of well-understood stochastic differential equations
(SDEs). This is best described by a metaphor from physics: take the
previous equation and assume $ \dim H = 1 $, $ \alpha(r) = 0 $ and $
\sigma(r) = \sigma $ a constant, i.e. an Ornstein-Uhlenbeck process
$$
dr_t = A r_t dt + \sigma dW_t, \, r_0 \in \mathbb{R}
$$
in dimension one describing the trajectory of a Brownian particle in
a (damping) velocity field $ x \mapsto A x $. If we move our
coordinate frame according to the vector field $ x \mapsto  A x $ we
observe a transformed movement of the particle, namely
$$
df_t = \exp(-At) \sigma dB_t, \, f_0 = r_0,
$$
which corresponds to a Brownian motion with time-dependent
volatility, since space is scaled by a factor $ \exp(-At) $ at time
$ t $ and the speed of the movement of the coordinate frame makes
the drift disappear. Loosely speaking, one ``jumps on the moving
frame'', where the speed of the frame is chosen equal to the drift.
In finite dimensions the advantage of this procedure is purely
conceptual, since analytically the both equations can be equally
well treated. If one imagines for a moment the same procedure for an
SPDE the advantage is much more than conceptual, since the
transformed equation, seen from the moving frame, is rather an SDE
than an SPDE, as the non-continuous drift term disappears in the
moving frame. More precisely, considering the variation of constants formula
\begin{align*}
r_t = S_t r_0 + \int_0^t S_{t-s} \alpha(r_s) ds + \int_0^t S_{t-s} \sigma(r_s)dW_s
\end{align*}
we recognize the dynamics of the transformed SDE, namely the process $f_t = S_{-t} r_t$ satisfies
\begin{align*}
df_t = S_{-t} \alpha(S_t f_t) dt + S_{-t} \sigma(S_t f_t) dW_t, \, f_0 = r_0.
\end{align*}

At this point it is clear that the drift term in infinite dimensions does usually not allow movements in negative
time direction, which is crucial for the approach. This limitation can be overcome by the Sz\H{o}kefalvi-Nagy theorem, which allows for
group extensions of given (pseudo-contractive) semigroups of linear operators. We emphasize that we do not need the particular structure
of this extension, which might be quite involved. The emphasis of this article is to provide a self-contained outline of this method in the realm
of jump-diffusions with path-dependent coefficients, which has not been treated in the literature so far.

Therefore we suggest the following approach to SPDEs, which is the guideline through this article:
\begin{itemize}
 \item consider the SDE obtained by transforming the SPDE with a time-dependent transformation $ r \mapsto S_{-t} r $ (jump to the moving frame).
 \item solve the transformed SDE.
 \item transform the solution process by $ r \mapsto S_t r $ in order to obtain a mild solution of the original SPDE (leave the moving frame).
\end{itemize}

In \cite{Ruediger-mild} and \cite{Marinelli-Prevot-Roeckner}
existence, uniqueness and regular dependence on initial data are
considered for SPDEs driven by a Wiener processes and Poisson random
measures. The authors also apply the Sz\H{o}kefalvi-Nagy theorem to
prove certain inequalities, which are crucial for their
considerations. In contrast our approach means that we reduce all
these separate considerations to the analysis of \emph{one}
transformed SDE, which corresponds then -- by means of the
time-dependent transformation -- to the solution of the given SPDE.

Our approach is based on the general jump-diffusion approach to
stochastic partial differential equations as presented in
\cite{Ruediger-mild} or \cite{Marinelli-Prevot-Roeckner}. In
contrast, our vector fields can be path-dependent in a general
sense, not only random as supposed in
\cite{Marinelli-Prevot-Roeckner}. Applications of this setting can
be found in recent work on volatility surfaces, where random
dependence of the vector fields is not enough. We first do the
obvious proofs for stochastic differential equations with values in
(separable) Hilbert spaces. Then we show that by our transformation
method (``jump to the moving frame'') we can transfer those results
to stochastic partial differential equations. In a completely
similar way we could have taken the setting for stochastic
differential equations in Ph.~Protter's book \cite{Protter}, which
is based on semi-martingales as driving processes and where we can
literally transfer the respective theorems into the setting of
stochastic partial differential equations. In particular all
$L^p$-estimates -- as extensively proved in \cite{Protter} -- can be
transferred into the setting of stochastic partial differential
equations.

The ``moving frame approach'' is a particular case of methods, where pull-backs with respect to flows are applied. Those methods have quite a long history in the theory of ODEs, PDEs and SDEs (pars pro toto we mention the Doss-Sussman method as described in \cite{RogWil:00} and the further material therein). In the realm of SPDEs the ``pull-back'' method has been successfully applied in \cite{BrzCapFla:88} with respect to noise vector fields. See also a discussion in \cite{BrzvNeVerWei:08} where this point of view is applied again, but a pull-back with respect to the PDE part has not been applied yet.

We shall now provide a guideline for the remainder of the article. In Section \ref{sec-integration} we define the fundamental concepts, notions and notations for stochastic integration with respect to Wiener processes and Poisson random measures. In Section \ref{sec-SDE-ex} and \ref{sec-SDE-lp} we provide for the sake of completeness existence and uniqueness results for Hilbert spaces valued SDEs and respective $ L^p $-estimates. In Sections \ref{sec-SDE-stability} and \ref{sec-SDE-regularity} we provide stability and regularity results for those SDEs. Section \ref{sec-SPDE-def} we introduce all necessary solution concepts for (semi-linear) SPDEs. In Section \ref{sec-SPDE-ex} we apply our method of the moving frame to existence and uniqueness questions. Section \ref{sec-stability-SPDE} is devoted to the study of stability and regularity for SPDEs. Section \ref{sec-SPDE-markov} and Section \ref{sec-SPDE-numerics} describe Markovian SPDE problems and several high order numerical schemes for SPDEs in this case. Again for the sake of completeness we provide a stochastic Fubini theorem with respect to compensated Poisson random measures in Appendix \ref{a}.

\section{Stochastic integration in Hilbert spaces}\label{sec-integration}

In this section, we shall outline the notion of stochastic integrals with
respect to an infinite dimensional Wiener process and with respect
to a compensated Poisson random measure. The construction of the
stochastic integral with respect to a Brownian motion follows
\cite[Sec. 4.2]{Da_Prato}. The construction of the stochastic
integral with respect to a Poisson measure is similar and can be
found in \cite{Barbara-Integration} or \cite[Sec. 2]{Knoche2}.

\subsection{Setting and Definitions}
From now on, let $(\Omega,\mathcal{F},(\mathcal{F}_t)_{t \geq
0},\mathbb{P})$ be a filtered probability space satisfying the usual
conditions. Furthermore, let $H$ denote a separable Hilbert space
with inner product $\langle \cdot,\cdot \rangle_H$ and associated
norm $\| \cdot \|_H$. If there is no ambiguity, we shall simply
write $\langle \cdot,\cdot \rangle$ and $\| \cdot \|$.

In the sequel, $\mathcal{P}$ denotes the predictable
$\sigma$-algebra on $\mathbb{R}_+$ and $\mathcal{P}_T$ denotes
predictable $\sigma$-algebra on $[0,T]$ for an arbitrary $T \in
\mathbb{R}_+$. We denote by $\lambda$ the Lebesgue measure on
$\mathbb{R}$.

For an arbitrary $p \geq 1$ and a finite time horizon $T \in
\mathbb{R}_+$ we define
\begin{align*}
L_T^p(\lambda;H) &:= L^p(\Omega \times [0,T], \mathcal{P}_T,
\mathbb{P} \otimes \lambda;H)
\end{align*}
and let $L^p(\lambda;H)$ be the space of all predictable process
$\Phi : \Omega \times \mathbb{R}_+ \rightarrow H$ such that for each
$T \in \mathbb{R}_+$ the restriction of $\Phi$ to $\Omega \times
[0,T]$ belongs to $L_T^p(\lambda;H)$. Furthermore, ${\mathcal{L}}_{\rm
loc}^p(\lambda;H)$ denotes the space of all predictable processes
$\Phi : \Omega \times \mathbb{R}_+ \rightarrow H$ such that
\begin{align*}
\mathbb{P} \bigg( \int_0^T \| \Phi_t \|^p dt < \infty \bigg) = 1
\quad \text{for all $T \in \mathbb{R}_+$.}
\end{align*}
Clearly, for each $\Phi \in \mathcal{L}_{\rm loc}^p(\lambda;H)$ the
path-by-path Stieltjes integral $\int_0^t \Phi_s ds$ exists.

Let $M_T^2(H)$ be the space of all square-integrable c\`adl\`ag
martingales $M : \Omega \times [0,T] \rightarrow H$, where
indistinguishable processes are identified. Endowed with the inner
product
\begin{align*}
(M,N) \mapsto \mathbb{E}[\langle M_T, N_T \rangle],
\end{align*}
the space $M_T^2(H)$ is a Hilbert space. The space $M_T^{2,c}(H)$,
consisting of all continuous elements from $M_T^2(H)$, is a closed
subspace of $M_T^2(H)$, which is a consequence of Doob's martingale
inequality \cite[Thm. 3.8]{Da_Prato}.

\subsection{Stochastic Integration with respect to Wiener
processes}\label{sec-integral-Wiener}

Let $U$ be another separable Hilbert space and $Q \in L(U)$ be a
compact, self-adjoint, strictly positive linear operator. Then there
exist an orthonormal basis $\{ e_j \}$ of $U$ and a bounded sequence
$\lambda_j$ of strictly positive real numbers such that
\begin{align*}
Qu = \sum_j \lambda_j \langle u,e_j \rangle e_j, \quad u \in U
\end{align*}
namely, the $\lambda_j$ are the eigenvalues of $Q$, and each $e_j$
is an eigenvector corresponding to $\lambda_j$, see, e.g.,
\cite[Thm. VI.3.2]{Werner}.

The space $U_0 := Q^{\frac{1}{2}}(U)$, equipped with inner product
$\langle u,v \rangle_{U_0} := \langle Q^{-\frac{1}{2}} u,
Q^{-\frac{1}{2}} v \rangle_U$, is another separable Hilbert space
and $\{ \sqrt{\lambda_j} e_j \}$ is an orthonormal basis.

Let $W$ be a $Q$-Wiener process \cite[p. 86,87]{Da_Prato}. We assume
that ${\rm tr}(Q) = \sum_j \lambda_j < \infty$. Otherwise, which
is the case if $W$ is a cylindrical Wiener process, there always
exists a separable Hilbert space $U_1 \supset U$ on which $W$ has a
realization as a finite trace class Wiener process, see \cite[Chap.
4.3]{Da_Prato}.

We denote by $L_2^0 := L_2(U_0,H)$ the space of Hilbert-Schmidt
operators from $U_0$ into $H$, which, endowed with the
Hilbert-Schmidt norm
\begin{align*}
\| \Phi \|_{L_2^0} := \sqrt{\sum_j \lambda_j \| \Phi e_j \|^2},
\quad \Phi \in L_2^0
\end{align*}
itself is a separable Hilbert space.

Following \cite[Chap. 4.2]{Da_Prato}, we define the stochastic
integral $\int_0^t \Phi_s dW_s$ as an isometry, extending the
obvious isometry on simple predictable processes, from
$L_T^2(W;L_2^0)$ to $M_T^{2,c}(H)$, where
\begin{align*}
L_T^2(W;L_2^0) := L^2(\Omega \times [0,T], \mathcal{P}_T, \mathbb{P}
\otimes \lambda; L_2^0).
\end{align*}
In particular, we obtain the \textit{It\^o-isometry}
\begin{align}\label{Ito-isometry-Wiener}
\mathbb{E} \left[ \bigg\| \int_0^t \Phi_s dW_s \bigg\|^2 \right] =
\mathbb{E} \bigg[ \int_0^t \| \Phi_s \|_{L_2^0}^2 ds \bigg], \quad t
\in [0,T]
\end{align}
for all $\Phi \in L_T^2(W;L_2^0)$. In a straightforward manner, we
extend the stochastic integral to the space $L^2(W;L_2^0)$ of all
predictable processes $\Phi : \Omega \times \mathbb{R}_+ \rightarrow
L_0^2$ such that the restriction of $\Phi$ to $\Omega \times [0,T]$
belongs to $L_T^2(W;L_2^0)$ for all $T \in \mathbb{R}_+$, and,
furthermore, to the space $\mathcal{L}_{\rm loc}^2(W;L_2^0)$ consisting of all
predictable processes $\Phi : \Omega \times \mathbb{R}_+ \rightarrow
L_0^2$ such that
\begin{align*}
\mathbb{P} \bigg( \int_0^T \| \Phi_t \|_{L_2^0}^2 dt < \infty \bigg)
= 1 \quad \text{for all $T \in \mathbb{R}_+$.}
\end{align*}
The integral process is unique up to indistinguishability.

There is an alternative view on the stochastic integral, which we
shall use in this text. According to \cite[Prop. 4.1]{Da_Prato}, the
sequence of stochastic processes $\{ \beta^j \}$ defined as $\beta^j :=
\frac{1}{\sqrt{\lambda_j}} \langle W, e_j \rangle$ is a sequence of
real-valued independent $(\mathcal{F}_t)$-Brownian motions and we
have the expansion
\begin{align*}
W = \sum_j \sqrt{\lambda _j} \beta^j e_j,
\end{align*}
where the series is convergent in the space $M^2(U)$ of $U$-valued
square-integrable martingales. Let $\Phi \in \mathcal{L}_{\rm loc}^2(W;L_2^0)$
be arbitrary. For each $j$ we set $\Phi^j := \sqrt{\lambda_j} \Phi
e_j$. Then we have
\begin{align*}
\int_0^t \Phi_s dW_s = \sum_j \int_0^t \Phi_s^j d\beta_s^j, \quad t
\in \mathbb{R}_+
\end{align*}
where the convergence is uniformly on compact time intervals in
probability, see \cite[Thm. 4.3]{Da_Prato}.

\subsection{Stochastic Integration with respect to Poisson random
measures}\label{sec-integral-Poisson}

Let $(E,\mathcal{E})$ be a measurable space which we assume to be a
\textit{Blackwell space} (see \cite{Dellacherie,Getoor}). We remark
that every Polish space with its Borel $\sigma$-field is a Blackwell
space.

Now let $\mu$ be a homogeneous Poisson random measure on $\mathbb{R}_+
\times E$, see \cite[Def. II.1.20]{JS}. Then its compensator is of
the form $dt \otimes F(dx)$, where $F$ is a $\sigma$-finite measure
on $(E,\mathcal{E})$.

We define the It\^o-integral $\int_0^t \int_E \Phi(s,x) (\mu(ds,dx)
- F(dx) ds)$ as an isometry, which extends the obvious isometry on
simple predictable processes, from $L_T^2(\mu;H)$ to $M_T^2(H)$,
where
\begin{align}\label{integrands-poisson_square_integrable}
L_T^2(\mu;H) := L^2(\Omega \times [0,T] \times E, \mathcal{P}_T
\otimes \mathcal{E}, \mathbb{P} \otimes \lambda \otimes F;H).
\end{align}
In particular, for each $\Phi \in L_T^2(\mu;H)$ we obtain the
\textit{It\^o-isometry}
\begin{align}\label{Ito-isometry-mu}
\mathbb{E} \left[ \bigg\| \int_0^t \int_E \Phi(s,x) (\mu(ds,dx) -
F(dx) ds) \bigg\|^2 \right] = \mathbb{E} \bigg[ \int_0^t \int_E \|
\Phi(s,x) \|^2 F(dx) ds \bigg]
\end{align}
for all $t \in [0,T]$. In a straightforward manner, we extend the
stochastic integral to the space $L^2(\mu;H)$ of all predictable
processes $\Phi : \Omega \times \mathbb{R}_+ \times E \rightarrow H$
such that the restriction of $\Phi$ to $\Omega \times [0,T] \times
E$ belongs to $L_T^2(\mu;H)$ for all $T \in \mathbb{R}_+$, and,
furthermore, to the space $\mathcal{L}_{\rm loc}^2(\mu;H)$ consisting of all
predictable processes $\Phi : \Omega \times \mathbb{R}_+ \times E \rightarrow H$ such that
\begin{align*}
\mathbb{P} \bigg( \int_0^T \| \Phi(t,x) \|^2 F(dx) dt < \infty
\bigg) = 1 \quad \text{for all $T \in \mathbb{R}_+$.}
\end{align*}
The integral process is unique up to indistinguishability.

Such a construction of the stochastic integral can, e.g., be found
in \cite[Sec. 4]{Applebaum} for the finite dimensional case and in
\cite{Barbara-Integration}, \cite[Sec. 2]{Knoche2} for the infinite
dimensional case.

\subsection{Path properties of stochastic
integrals}\label{sec-path-properties}

It is apparent that for every $\Phi \in \mathcal{L}_{\rm loc}^p(\lambda;H)$,
where $p \geq 1$, the path-by-path Stieltjes integral
$\int_0^{\bullet} \Phi_s ds$ has continuous sample paths.

As outlined in Section \ref{sec-integral-Wiener}, we have first
defined the stochastic integral $\int_0^t \Phi_s dW_s$ as an
isometry from $L_T^2(W;L_2^0)$ to $M_T^{2,c}(H)$, the space of all
square-integrable continuous martingales, and then extended it by
localization. Therefore, for each $\Phi \in \mathcal{L}_{\rm loc}^2(W;L_2^0)$,
the trajectories of the integral process $\int_0^{\bullet} \Phi_s
dW_s$ are continuous.

Similarly, the stochastic integral $\int_0^t \int_E \Phi(s,x)
(\mu(ds,dx) - F(dx) ds)$, outlined in Section
\ref{sec-integral-Poisson}, is, in the first step, defined as an
isometry from $L_T^2(\mu;H)$ to $M_T^2(H)$, the space of all
square-integrable c\`adl\`ag martingales, and then extended by
localization. Hence, for each $\Phi \in \mathcal{L}_{\rm loc}^2(\mu;H)$ the
integral process $\int_0^{\bullet} \int_E \Phi(s,x) (\mu(ds,dx) -
F(dx) ds)$ has c\`adl\`ag sample paths.

\subsection{Independence of the driving terms}\label{sec-independence}

We remark that the Wiener process $W$ and the Poisson random measure $\mu$ are independent, which we will actually only need in Section \ref{sec-SPDE-numerics}.

The asserted independence is provided by using the semimartingale theory from Jacod and Shiryaev \cite{JS}. Indeed, for a continuous local martingale $M$ and a purely discontinuous local martingale $N$, which are both assumed to be processes with independent increments and both considered with respect to the same filtration, the semimartingale $X = (M,N)$ is again a process with independent increments, because its semimartingale characteristics (see \cite[Def. II.2.6]{JS}), which we can easily compute from those of $M$ and $N$, are also deterministic. Here we need the fact that $ (M,N) = (M,0) + (0,N) $ is a decomposition into a continuous and purely discontinuous local martingale. Computing the characteristic functions of $M$, $N$ and $X$ by means of \cite[Thm. II.4.15]{JS} yields the desired independence.

\section{Existence and uniqueness of solutions for stochastic differential equations}\label{sec-SDE-ex}

Since we shall show that -- in case of pseudo-contractive strongly continuous semigroups -- it is equivalent to consider SPDEs on the one hand
or time-dependent SDEs on the other hand, we need the basic results for time-dependent SDEs with possibly infinite dimensional state space
at hand. In this section we prove existence and uniqueness results for stochastic differential equations (SDEs) on a possibly
infinite dimensional state space. 
The results are fairly standard, but we provide them in order to keep our presentation
self-contained and to introduce certain notation which we shall need in the further sections.

For an interval $I \subset \mathbb{R}_+$ we define the space $C(I;H) := C(I;L^2(\Omega;H))$ of
all continuous functions from $I$ into $L^2(\Omega;H)$. If the interval $I$ is compact, then $C(I;H)$
is a Banach space with respect to the norm
\begin{align*}
\| r \|_{I} := \sup_{t \in I} \| r_t \|_{L^2(\Omega;H)} = \sqrt{
\sup_{t \in I} \mathbb{E} [ \| r_t \|^2 ] }.
\end{align*}
Note that $C(I;H)$ is a space consisting of continuous curves of
equivalence classes of random variables. For each element $r \in C(I;H)$ we can associate an $H$-valued, mean-square continuous
process $\tilde{r} = (\tilde{r}_t)_{t \in I}$, which is unique up to a version.

Let $C_{\rm ad}(I;H)$ be the subspace consisting of all adapted curves
from $C(I;H)$. Note that, by the
completeness of the filtration $(\mathcal{F}_t)_{t \geq 0}$,
adaptedness of a curve $r \in C(I;H)$ is independent of the choice
of the representative. If the interval $I$ is compact, then the subspace $C_{\rm ad}(I;H)$ is closed with respect to the norm $\| \cdot \|_I$.

We shall also consider the spaces $\mathcal{C}(I;H)$ and $\mathcal{C}_{\rm ad}(I;H)$ of all mean-square continuous and of all adapted, mean-square continuous processes $r \in C(I;\mathcal{L}^2(\Omega;H))$. Note that for each $r \in \mathcal{C}(I;H)$ the equivalence class $[r]$ belongs to $C(I;H)$, and if
$r \in \mathcal{C}_{\rm ad}(I;H)$, then we have $[r] \in C_{\rm ad}(I;H)$.

If no confusion concerning the Hilbert space $H$ is possible, we shall use the abbreviations $C(I)$, $C_{\rm ad}(I)$, $\mathcal{C}(I)$ and $\mathcal{C}_{\rm ad}(I)$ for $C(I;H)$, $C_{\rm ad}(I;H)$, $\mathcal{C}(I;H)$ and $\mathcal{C}_{\rm ad}(I;H)$.

We denote by $H_{\mathcal{P}}$ resp. $H_{\mathcal{P} \otimes
\mathcal{E}}$ the space of all predictable processes $r : \Omega
\times \mathbb{R}_+ \rightarrow H$ resp. $r : \Omega \times
\mathbb{R}_+ \times E \rightarrow H$.

We shall now deal with stochastic differential equations of the kind
\begin{align}\label{equation-strong}
\left\{
\begin{array}{rcl}
dr_t & = & \alpha(r)_t dt + \sigma(r)_t dW_t + \int_E \gamma(r)(t,x)
(\mu(dt,dx) - F(dx) dt) \medskip
\\ r|_{[0,t_0]} & = & h,
\end{array}
\right.
\end{align}
where $\alpha : C_{\rm ad}(\mathbb{R}_+) \rightarrow H_{\mathcal{P}}$, $\sigma
: C_{\rm ad}(\mathbb{R}_+) \rightarrow (L_2^0)_{\mathcal{P}}$ and $\gamma :
C_{\rm ad}(\mathbb{R}_+) \rightarrow H_{\mathcal{P} \otimes \mathcal{E}}$. Fix $t_0 \in \mathbb{R}_+$ and $h \in \mathcal{C}_{\rm ad}[0,t_0]$.

\begin{definition}\label{def-solution-SDE}
A process $r \in \mathcal{C}_{\rm ad}(\mathbb{R}_+)$ is called a {\rm solution} for
(\ref{equation-strong}) if we have $r|_{[0,t_0]} = h$, $\alpha([r]) \mathbbm{1}_{[t_0,\infty)} \in
\mathcal{L}_{\rm loc}^1(\lambda;H)$, $\sigma([r]) \mathbbm{1}_{[t_0,\infty)} \in \mathcal{L}_{\rm loc}^2(W;L_0^2)$,
$\gamma([r]) \mathbbm{1}_{[t_0,\infty)} \in \mathcal{L}_{\rm loc}^2(\mu;H)$ and almost surely
\begin{equation}\label{SDE-solved-eqn}
\begin{aligned}
r_t &= h_{t_0} + \int_{t_0}^t \alpha([r])_s ds + \int_{t_0}^t \sigma([r])_s dW_s
\\ &\quad + \int_{t_0}^t \int_E \gamma([r])(s,x) (\mu(ds,dx) - F(dx) ds), \quad t \geq
t_0.
\end{aligned}
\end{equation}
As pointed out in Section \ref{sec-integration}, the stochastic integrals at the right-hand side of (\ref{SDE-solved-eqn}) are only determined up to indistinguishability. Therefore, {\em uniqueness} of solutions for (\ref{equation-strong}) is meant up to indistinguishability on the interval $[t_0,\infty)$, that
is, for two solutions $r, \tilde{r} \in \mathcal{C}_{\rm ad}(\mathbb{R}_+)$ we have $\mathbb{P}(\bigcap_{t \geq t_0} \{ r_t = \tilde{r}_t \}) = 1$.
\end{definition}

\begin{remark}
Note that in this definition time-dependence of the vector fields is
naturally included into the setting. Also observe that for $t_0 = 0$ we have $C_{\rm ad}[0,t_0] = L^2(\Omega,\mathcal{F}_0,\mathbb{P};H)$ and $\mathcal{C}_{\rm ad}[0,t_0] = \mathcal{L}^2(\Omega,\mathcal{F}_0,\mathbb{P};H)$.
\end{remark}

The following standard assumptions are crucial for existence and uniqueness:

\begin{assumption}\label{ass-0-strong}
We assume that for all $T \in \mathbb{R}_+$ and all $r^1,r^2 \in C_{\rm ad}(\mathbb{R}_+)$ with $r^1|_{[0,T]} = r^2|_{[0,T]}$ we have
\begin{align*}
\alpha(r^1)|_{[0,T]} &= \alpha(r^2)|_{[0,T]},
\\ \sigma(r^1)|_{[0,T]} &= \sigma(r^2)|_{[0,T]},
\\ \gamma(r^1)|_{[0,T] \times E} &= \gamma(r^2)|_{[0,T] \times E}.
\end{align*}
\end{assumption}

\begin{assumption}\label{ass-1-strong}
Denoting by $\mathbf{0} \in C_{\rm ad}(\mathbb{R}_+)$ the zero process, we
assume that
\begin{align}\label{ass-alpha-regular}
t \mapsto \mathbb{E} [ \| \alpha({\mathbf{0}})_t \|^2 ] &\in \mathcal{L}_{\rm
loc}^{1}(\mathbb{R}_+),
\\ \label{ass-sigma-regular} t \mapsto \mathbb{E} [ \| \sigma({\mathbf{0}})_t \|_{L_2^0}^2 ]
&\in \mathcal{L}_{\rm loc}^{1}(\mathbb{R}_+),
\\ \label{ass-gamma-regular} t \mapsto \mathbb{E} \bigg[ \int_E \| \gamma({\mathbf{0}})(t,x)
\|^2 F(dx) \bigg] &\in \mathcal{L}_{\rm loc}^{1}(\mathbb{R}_+).
\end{align}
\end{assumption}

\begin{assumption}\label{ass-2-strong}
We assume there is a function
\begin{align}\label{L-function-integrable}
L \in \mathcal{L}_{\rm loc}^{2}(\mathbb{R}_+)
\end{align}
such that for all $t \in \mathbb{R}_+$ we have
\begin{align}\label{Lip-alpha-SDE-strong}
\mathbb{E} [ \| \alpha(r^1)_t - \alpha(r^2)_t \|^2 ] &\leq L(t)^2 \| r^1 - r^2 \|_{[0,t]}^2,
\\ \label{Lip-sigma-SDE-strong} \mathbb{E} [ \|
\sigma(r^1)_t - \sigma(r^2)_t \|_{L_2^0}^2 ] &\leq L(t)^2 \| r^1 - r^2 \|_{[0,t]}^2,
\\ \label{Lip-gamma-SDE-strong} \mathbb{E} \bigg[ \int_E \| \gamma(r^1)(t,x)
- \gamma(r^2)(t,x) \|^2 F(dx) \bigg] &\leq L(t)^2 \| r^1 - r^2 \|_{[0,t]}^2
\end{align}
for all $r^1,r^2 \in C_{\rm ad}(\mathbb{R}_+)$.
\end{assumption}

\begin{remark}
For $p \geq 1$ the space $\mathcal{L}_{\rm loc}^{p}(\mathbb{R}_+)$ denotes the
space of all measurable functions $f : \mathbb{R}_+ \rightarrow
\mathbb{R}$ such that the restriction $f|_{[0,T]}$ belongs to
$\mathcal{L}^p[0,T]$ for every $T \in \mathbb{R}_+$. Note that
(\ref{ass-alpha-regular}), (\ref{ass-sigma-regular}),
(\ref{ass-gamma-regular}) and (\ref{L-function-integrable}) are in
particular satisfied if the respective functions are bounded on
compact intervals.
\end{remark}

\begin{remark}
Note that Assumptions \ref{ass-0-strong}, \ref{ass-1-strong}, \ref{ass-2-strong} are satisfied for a wide class of SDEs (and thus -- by the method of the moving frame -- SPDEs) with path-dependent coefficients. As an example, we will consider equations with characteristic coefficients depending on the randomness $\omega$, the time $t$ and finitely many states from the path on the interval $[0,t]$, see Corollary \ref{cor-several} below, which together with Remark \ref{remark-Lp-state} (see also Remark \ref{remark-dissipative}) generalizes \cite[Thm. 2.4]{Marinelli-Prevot-Roeckner}. We also emphasize that the Lipschitz function $L$ only needs to be locally square-integrable.
\end{remark}

\begin{remark}
In the book of Ph. Protter \cite{Protter} stochastic differential equations driven by semimartingales are studied. The characteristic coefficients are mappings $F : \mathbb{D} \rightarrow \mathbb{D}$, where $\mathbb{D}$ denotes the space of adapted c\`{a}dl\`{a}g processes. In \cite[Thm. V.7]{Protter} they are assumed to be {\em functional Lipschitz}, i.e. for any $X,Y \in \mathbb{D}$ we have
\begin{align}\label{functional-tau}
F(X)^{\tau-} = F(Y)^{\tau-} \quad \text{for any stopping time $\tau$ with $X^{\tau-} = Y^{\tau-}$}
\end{align}
and almost surely
\begin{align}\label{functional-Lipschitz}
\| F(X)_t - F(Y)_t \| \leq K_t \sup_{s \in [0,t]} \| X_s - Y_s \| \quad \text{for each $t \geq 0$,}
\end{align}
where $K = (K_t)_{t \geq 0}$ is an increasing (finite) process. By localization, Protter \cite{Protter} assumes that $K$ is uniformly bounded by some finite constant $k > 0$, see \cite[Lemmas V.1, V.2]{Protter}. Taking expectation in (\ref{functional-Lipschitz}) then yields
\begin{align}\label{Lip-topology}
\mathbb{E} [ \| F(X)_t - F(Y)_t \|^2 ] \leq k^2 \mathbb{E} \bigg[ \sup_{s \in [0,t]} \| X_s - Y_s \|^2 \bigg] = k^2 \| X-Y \|_{S^2[0,t]}^2,
\end{align}
and existence and uniqueness is proven by a fixed point argument on the space $S^2$. We, in contrast, will apply a fixed point argument on the space $C_{\rm ad}(\mathbb{R}_+)$, and show the existence of a c\`{a}dl\`{a}g version afterwards, see Theorem \ref{thm-existence} below. Note that Assumption \ref{ass-0-strong} corresponds to (\ref{functional-tau}) and Assumption \ref{ass-2-strong} corresponds to (\ref{Lip-topology}). Hence, our assumptions can be regarded as an analogue to the functional Lipschitz property in \cite{Protter}.
\end{remark}

\begin{lemma}\label{lemma-well-defined}
For each $r \in C_{\rm ad}(\mathbb{R}_+)$ the functions
\begin{align}\label{alpha-cont}
&t \mapsto \mathbb{E} \bigg[ \int_{0}^{t} \| \alpha(r)_s \|^2 ds
\bigg],
\\ \label{sigma-cont} &t \mapsto \mathbb{E} \bigg[ \int_{0}^{t} \| \sigma(r)_s \|_{L_2^0}^2 ds \bigg],
\\ \label{gamma-cont} &t \mapsto \mathbb{E} \bigg[ \int_{0}^{t} \int_E \| \gamma(r)(s,x) \|^2 F(dx) ds
\bigg]
\end{align}
are well-defined and continuous on $\mathbb{R}_+$.
\end{lemma}

\begin{proof}
Let $r \in C_{\rm ad}(\mathbb{R}_+)$ and $t \in \mathbb{R}_+$ be arbitrary.
Using the Lipschitz conditions (\ref{Lip-alpha-SDE-strong}),
(\ref{Lip-sigma-SDE-strong}), (\ref{Lip-gamma-SDE-strong}) we obtain
\begin{align*}
&\mathbb{E} \bigg[ \int_{0}^{t} \| \alpha(r)_s \|^2 ds \bigg] \leq 2
\int_{0}^{t} L(s)^2 \| r \|_{[0,s]}^2 ds + 2 \int_{0}^{t}
\mathbb{E}[\| \alpha(\mathbf{0})_s \|^2]ds,
\\ &\mathbb{E} \bigg[ \int_{0}^{t} \| \sigma(r)_s \|_{L_2^0}^2 ds \bigg]
\leq 2 \int_{0}^{t} L(s)^2 \| r \|_{[0,s]}^2 ds + 2
\int_{0}^{t} \mathbb{E}[\| \sigma(\mathbf{0})_s \|_{L_2^0}^2]ds,
\\ &\mathbb{E} \bigg[ \int_{0}^{t} \int_E \| \gamma(r)(s,x) \|^2
F(dx) ds \bigg]
\\ &\leq 2 \int_{0}^{t} L(s)^2 \| r \|_{[0,s]}^2 ds + 2
\int_{0}^{t} \mathbb{E}\bigg[ \int_E \| \gamma(\mathbf{0})(s,x) \|^2
F(dx) \bigg] ds.
\end{align*}
Note that, by (\ref{L-function-integrable}), we have
\begin{align*}
\int_{0}^{t} L(s)^2 \| r \|_{[0,s]}^2 ds \leq \| r \|_{[0,t]}^2
\int_{0}^{t} L(s)^2 ds < \infty.
\end{align*}
Together with (\ref{ass-alpha-regular}), (\ref{ass-sigma-regular}),
(\ref{ass-gamma-regular}) we deduce that the functions in
(\ref{alpha-cont}), (\ref{sigma-cont}), (\ref{gamma-cont}) are
well-defined. The continuity follows from Lebesgue's theorem.
\end{proof}

According to Lemma \ref{lemma-well-defined}, for all $r \in C_{\rm ad}(\mathbb{R}_+)$ we have $\alpha(r) \in
L^2(\lambda;H)$, $\sigma(r) \in L^2(W;L_2^0)$ and $\gamma(r) \in
L^2(\mu;H)$. This ensures that the following stochastic integrals in (\ref{def-Lambda}) are well-defined.

For any $T \in \mathbb{R}_+$ and $r \in C_{\rm ad}[0,T]$ we define 
\begin{align*}
\alpha(r) &:= \alpha(\tilde{r})|_{[0,T]},
\\ \sigma(r) &:= \sigma(\tilde{r})|_{[0,T]},
\\ \gamma(r) &:= \gamma(\tilde{r})|_{[0,T] \times E},
\end{align*}
where we have chosen $\tilde{r} \in C_{\rm ad}(\mathbb{R}_+)$ such that $r = \tilde{r}|_{[0,T]}$. Such an element $\tilde{r}$ always exists. Take, for example, the constant continuation $\tilde{r}_t := r_T$ for $t \geq T$. Notice also that this definition is independent of the choice of $\tilde{r}$ by virtue of Assumption \ref{ass-0-strong}.

Let us fix $t_0 \in \mathbb{R}_+$, $T \geq t_0$, $h \in C_{\rm ad}[0,t_0]$ and $r \in
C_{\rm ad}[0,T]$. We define $\Lambda_{h}(r)$ by $\Lambda_{h}(r)|_{[0,t_0]} := h$ and
\begin{equation}\label{def-Lambda}
\begin{aligned}
\Lambda_{h}(r)_t &:= h_{t_0} + \int_{t_0}^t \alpha(r)_s ds + \int_{t_0}^t
\sigma(r)_s dW_s
\\ &\quad + \int_{t_0}^t \int_E \gamma(r)(s,x) (\mu(ds,dx) - F(dx) ds), \quad t \in [t_0,T].
\end{aligned}
\end{equation}
By H\"older's inequality, the
It\^o-isometries (\ref{Ito-isometry-Wiener}),
(\ref{Ito-isometry-mu}) and Lemma \ref{lemma-well-defined}, the process
$\Lambda_{h}(r)$ is mean-square continuous. By taking the respective equivalence classes,
this induces a mapping $\Lambda_{h} : C_{\rm ad}[0,T] \rightarrow C_{\rm ad}[0,T]$.

In an analogous fashion, we define a mapping $\Lambda_{h} : C_{\rm ad}(\mathbb{R}_+) \rightarrow C_{\rm ad}(\mathbb{R}_+)$.

Now we fix $T_1,T_2 \in \mathbb{R}_+$ with $T_1 \leq T_2$ and $r^1 \in C_{\rm ad}[0,T_1]$. For $r^2 \in C_{\rm ad}[T_1,T_2]$ we have 
\begin{align*}
(r^1,r^2) := \big( (r_s^1)_{s \in [0,T_1]},(r_s^2 + r_{T_1}^1 - r_{T_1}^2)_{s \in (T_1,T_2]} \big) \in C_{\rm ad}[0,T_2].
\end{align*}
Hence, we can define
\begin{align*}
\Gamma_{r^1}(r^2)_t &:= r_{T_1}^1 + \int_{T_1}^t \alpha ( r^1,r^2 )_s ds + \int_{T_1}^t \sigma ( r^1,r^2 )_s dW_s 
\\ &\quad + \int_{T_1}^t \gamma ( r^1,r^2 )(s,x) (\mu(ds,dx) - F(dx)ds), \quad t \in [T_1,T_2].
\end{align*}
By H\"older's inequality, the
It\^o-isometries (\ref{Ito-isometry-Wiener}),
(\ref{Ito-isometry-mu}) and Lemma \ref{lemma-well-defined}, the process
$\Gamma_{r^1}(r^2)$ is mean-square continuous. By taking the respective equivalence classes,
this induces a mapping $\Gamma_{r^1} : C_{\rm ad}[T_1,T_2] \rightarrow C_{\rm ad}[T_1,T_2]$.

\begin{lemma}\label{lemma-contraction-seq}
Let $t_0 \in \mathbb{R}_+$ be arbitrary.
There exists a sequence $t_0 = T_0 < T_1 < T_2 < \ldots$ with $T_n \rightarrow \infty$ such that for all $n \in \mathbb{N}_0$ and all $h \in C_{\rm ad}[0,T_n]$ the map $\Gamma_{h}$ is a contraction on $C_{\rm ad}[T_n,T_{n+1}]$.
\end{lemma}

\begin{proof}
We choose an arbitrary $\epsilon \in (0,1)$. Let $\delta > 0$ be such that
\begin{align}\label{f-small}
f(t) \leq \epsilon, \quad t \in [0,\delta]
\end{align}
where $f(t) := 12(t+2)$. By (\ref{L-function-integrable}) and Lebesgue's theorem, the map $g : \mathbb{R}_+ \rightarrow \mathbb{R}_+$, $g(t) = \int_0^t L(s)^2 ds$ is continuous.
Since $g$ is uniformly continuous on compact intervals of $\mathbb{R}_+$, there exists a sequence $t_0 = T_0 < T_1 < T_2 < \ldots$ with $\sup_{n \in \mathbb{N}_0} |T_{n+1} - T_n| \leq \delta$ and $T_n \rightarrow \infty$ such that
\begin{align}\label{g-small}
|g(T_n) - g(T_{n+1})| \leq \epsilon \quad \text{for all $n \in \mathbb{N}_0$.}
\end{align}
Let $n \in \mathbb{N}_0$ and $h \in C_{\rm ad}[0,T_n]$ be arbitrary. We fix $r^1,r^2 \in C_{\rm ad}[T_n,T_{n+1}]$ and $t \in [T_n,T_{n+1}]$.
By using H\"older's inequality and (\ref{Lip-alpha-SDE-strong}) we obtain
\begin{align*}
&\mathbb{E} \Bigg[ \bigg\| \int_{T_n}^t ( \alpha(h,r^1)_s - \alpha(h,r^2)_s ) ds \bigg\|^2 \Bigg] 
\\ &\leq (t-T_n) \int_{T_n}^t L(s)^2 \| r^1 + h_{T_n} - r_{T_n}^1 - (r^2 + h_{T_n} - r_{T_n}^2 ) \|_{[T_n,s]}^2 ds
\\ &\leq 4 (t-T_n) \bigg( \int_{T_n}^t L(s)^2 ds \bigg) \| r^1 - r^2 \|_{[T_n,T_{n+1}]}^2.
\end{align*}
The It\^o-isometry (\ref{Ito-isometry-Wiener}) and
(\ref{Lip-sigma-SDE-strong}) yield
\begin{align*}
&\mathbb{E} \Bigg[ \bigg\| \int_{T_n}^t ( \sigma(h,r^1)_s - \sigma(h,r^2)_s ) dW_s \bigg\|^2 \Bigg] 
\\ &\leq \int_{T_n}^t L(s)^2 \| r^1 + h_{T_n} - r_{T_n}^1 - (r^2 + h_{T_n} - r_{T_n}^2 ) \|_{[T_n,s]}^2 ds
\\ &\leq 4 \bigg( \int_{T_n}^t L(s)^2 ds \bigg) \| r^1 - r^2 \|_{[T_n,T_{n+1}]}^2,
\end{align*}
and the It\^o-isometry (\ref{Ito-isometry-mu}) and
(\ref{Lip-gamma-SDE-strong}) give us an analogous estimate for the jump part.
Thus, we obtain for all $t \in [T_n,T_{n+1}]$ the estimate
\begin{align*}
\mathbb{E}[\| \Gamma_{h}(r^1)_t - \Gamma_{h}(r^2)_t \|^2] &\leq 12(t-T_n+2) \bigg( \int_{T_n}^t L(s)^2 ds \bigg) \| r^1 - r^2 \|_{[T_n,T_{n+1}]}^2
\\ &= f(t-t_n) (g(t) - g(T_n)) \| r^1 - r^2 \|_{[T_n,T_{n+1}]}^2,
\end{align*}
which implies, by taking into account (\ref{f-small}) and (\ref{g-small}),
\begin{align*}
\| \Gamma_{h}(r^1) - \Gamma_{h}(r^2) \|_{[T_n,T_{n+1}]} \leq \epsilon \| r^1 - r^2 \|_{[T_n,T_{n+1}]},
\end{align*}
proving that $\Gamma_{h}$ is a contraction on $C_{\rm ad}[T_n,T_{n+1}]$.
\end{proof}

\begin{theorem}\label{thm-existence}
Suppose that Assumptions \ref{ass-0-strong}, \ref{ass-1-strong}, \ref{ass-2-strong}
are fulfilled. Then, for each $t_0 \in \mathbb{R}_+$ and $h \in
\mathcal{C}_{\rm ad}[0,t_0]$ there exists a unique
solution $r \in \mathcal{C}_{\rm ad}(\mathbb{R}_+)$ for
(\ref{equation-strong}) with c\`{a}dl\`{a}g paths on $[t_0,\infty)$, and it satisfies
\begin{align}\label{solution-in-S2}
\mathbb{E} \bigg[ \sup_{t \in [t_0,T]} \| r_t \|^2 \bigg] < \infty
\quad \text{for all $T \geq t_0$.}
\end{align}
\end{theorem}

\begin{proof}
Let $t_0 \in \mathbb{R}_+$ and $h \in
\mathcal{C}_{\rm ad}[0,t_0]$ be arbitrary. We identify $h$ with its equivalence class and fix a sequence $(T_n)_{n \in \mathbb{N}}$ as in Lemma \ref{lemma-contraction-seq}. By induction we shall prove that for each $n \in \mathbb{N}_0$ the fixed point equation
\begin{align}\label{fixed-point-Lambda}
r^{n} = \Lambda_{h}(r^{n}), \quad r^{n} \in C_{\rm ad}[0,T_{n}]
\end{align}
has a unique solution. For $n=0$ the unique solution for (\ref{fixed-point-Lambda}) is given by $r^0 = h$.
We proceed with the induction step $n \rightarrow n+1$.
By the Banach fixed point theorem there exists a unique solution for
\begin{align}\label{fix-ind-step}
\tilde{r}^{n+1} = \Gamma_{r^n}(\tilde{r}^{n+1}), \quad \tilde{r}^{n+1} \in C_{\rm ad}[T_n,T_{n+1}].
\end{align}
The process $r^{n+1} := ((r^n)_{t \in [0,T_n]},(\tilde{r}^{n+1})_{t \in (T_n,T_{n+1}]})$ belongs to $C_{\rm ad}[0,T_{n+1}]$, because $r_{T_n}^n = \tilde{r}_{T_n}^{n+1}$ by (\ref{fix-ind-step}), and, by taking into account Assumption \ref{ass-0-strong}, it is the unique solution for
\begin{align*}
r^{n+1} = \Lambda_{h}(r^{n+1}), \quad r^{n+1} \in C_{\rm ad}[0,T_{n+1}].
\end{align*}
Since $T_n \rightarrow \infty$, there exists, by noting Assumption \ref{ass-0-strong} again, a unique solution $r \in C_{\rm ad}(\mathbb{R}_+)$
for the fixed point equation
\begin{align}\label{fixed-point-version}
r = \Lambda_{h}(r), \quad r \in C_{\rm ad}(\mathbb{R}_+).
\end{align}
The right-hand side of (\ref{fixed-point-version}) consists of the sum of stochastic integrals. Therefore, there exists a representative $\tilde{r} \in \mathcal{C}_{\rm ad}(\mathbb{R}_+)$ of $\Lambda_{h}(r)$ with c\`adl\`ag paths on $[t_0,\infty)$, see Section \ref{sec-path-properties}. Equation (\ref{fixed-point-version}) yields, up to indistinguishability,
\begin{align*}
\tilde{r}_t &= h_{t_0} + \int_{t_0}^t \alpha([\tilde{r}])_s ds + \int_{t_0}^t \sigma([\tilde{r}])_s dW_s
\\ &\quad + \int_{t_0}^t \int_E \gamma([\tilde{r}])(s,x) (\mu(ds,dx) - F(dx) ds), \quad t \geq t_0.
\end{align*}
Since any two representatives of $r$, which are c\`adl\`ag on $[t_0,\infty)$, are indistinguishable on $[t_0,\infty)$,
this shows that $\tilde{r}$ is the unique solution for
(\ref{equation-strong}). Relation (\ref{solution-in-S2})
is established by H\"older's inequality, Doob's martingale
inequality \cite[Thm. 3.8]{Da_Prato}, the It\^o-isometries
(\ref{Ito-isometry-Wiener}), (\ref{Ito-isometry-mu}) and Lemma
\ref{lemma-well-defined}.
\end{proof}

\begin{remark}
The idea work on the space $C_{\rm ad}(\mathbb{R}_+)$ already appears in the
proof of \cite[Thm. 4.1]{Onno}, which deals with infinite
dimensional stochastic differential equations driven by Wiener
processes.
\end{remark}

\section{$L^p$-estimates}\label{sec-SDE-lp}
In order to carry $L^p$-theory from SDEs with possibly infinite dimensional state space to SPDEs we provide the relevant results for SDEs
here. For the SDEs of Section \ref{sec-SDE-ex} the full theory of $L^p$-estimates for
solutions of stochastic differential equations holds true.

Let $p \geq 2$ be arbitrary. In this section, for any interval $I \subset \mathbb{R}_+$ we consider the space $C(I) := C(I;L^p(\Omega;H))$ of
all continuous functions from $I$ into $L^p(\Omega;H)$. If the interval $I$ is compact, we
equip $C(I)$ with the norm
\begin{align*}
\| r \|_{I} := \sup_{t \in I} \| r_t \|_{L^p(\Omega;H)} = \bigg(
\sup_{t \in I} \mathbb{E} [ \| r_t \|^p ] \bigg)^{\frac{1}{p}}.
\end{align*}
We replace Assumptions \ref{ass-1-strong} and \ref{ass-2-strong} by the following stronger assumptions.

\begin{assumption}\label{ass-1-Lp}
Denoting by $\mathbf{0} \in C_{\rm ad}(\mathbb{R}_+)$ the zero process, we
assume that
\begin{align*}
t \mapsto \mathbb{E} [ \| \alpha({\mathbf{0}})_t \|^p ] &\in \mathcal{L}_{\rm
loc}^{1}(\mathbb{R}_+),
\\ t \mapsto \mathbb{E} [ \| \sigma({\mathbf{0}})_t \|_{L_2^0}^p ] &\in \mathcal{L}_{\rm
loc}^{1}(\mathbb{R}_+),
\\ t \mapsto \mathbb{E} \bigg[ \int_E \|
\gamma({\mathbf{0}})(t,x) \|^p F(dx) \bigg] &\in \mathcal{L}_{\rm
loc}^{1}(\mathbb{R}_+).
\end{align*}
\end{assumption}

\begin{assumption}\label{ass-2-Lp}
We assume there is a function
\begin{align*}
L \in \mathcal{L}_{\rm loc}^{p}(\mathbb{R}_+)
\end{align*}
such that for all $t \in \mathbb{R}_+$ we have
\begin{align}\label{Lip-alpha-Lp}
\mathbb{E} [ \| \alpha(r^1)_t - \alpha(r^2)_t \|^p ] &\leq L(t)^p \| r^1 - r^2 \|_{[0,t]}^p,
\\ \label{Lip-sigma-Lp} \mathbb{E} [ \|
\sigma(r^1)_t - \sigma(r^2)_t \|_{L_2^0}^p ] &\leq L(t)^p \| r^1 - r^2 \|_{[0,t]}^p,
\\ \label{Lip-gamma-Lp} \mathbb{E} \bigg[ \int_E \| \gamma(r^1)(t,x)
- \gamma(r^2)(t,x) \|^p F(dx) \bigg]&
\\ \notag + \mathbb{E} \Bigg[ \bigg( \int_E \| \gamma(r^1)(t,x)
- \gamma(r^2)(t,x) \|^2 F(dx) \bigg)^{\frac{p}{2}} \Bigg] &\leq L(t)^p \|
r^1 - r^2 \|_{[0,t]}^p
\end{align}
for all $r^1,r^2 \in C_{\rm ad}(\mathbb{R}_+)$.
\end{assumption}

\begin{theorem}\label{thm-Lp-solutions}
Suppose that Assumptions \ref{ass-0-strong}, \ref{ass-1-Lp}, \ref{ass-2-Lp}
are fulfilled. Then, for each $t_0 \in \mathbb{R}_+$ and $h \in
\mathcal{C}_{\rm ad}[0,t_0]$ there exists a unique
solution $r \in \mathcal{C}_{\rm ad}(\mathbb{R}_+)$ for
(\ref{equation-strong}) with c\`adl\`ag paths on $[t_0,\infty)$, and it satisfies
\begin{align*}
\mathbb{E} \bigg[ \sup_{t \in [t_0,T]} \| r_t \|^p \bigg] < \infty
\quad \text{for all $T \geq t_0$.}
\end{align*}
\end{theorem}

The proof is established by applying the reasonings from the previous section directly.
We do not go into detail here, but indicate how we apply the Banach fixed point
theorem in this situation, which relies on Burkholder-Davis-Gundy and Bichteler-Jacod type arguments.

By using H\"older's inequality and
(\ref{Lip-alpha-Lp}) we obtain
\begin{align*}
&\mathbb{E} \left[ \bigg\| \int_{T_n}^t (\alpha(h,r^1)_s - \alpha(h,r^2)_s)
ds \bigg\|^p \right] 
\\ &\leq (t - T_n)^{p-1} \mathbb{E} \bigg[ \int_{T_n}^t \|
\alpha(h,r^1)_s - \alpha(h,r^2)_s \|^p ds \bigg]
\\ &\leq 2^p (t-T_n)^{p-1} \bigg( \int_0^t L(s)^p ds \bigg) \| r^1 - r^2 \|_{[T_n,T_{n+1}]}^p.
\end{align*}
By the Burkholder-Davis-Gundy inequality, H\"older's inequality
and (\ref{Lip-sigma-Lp}) we have
\begin{align*}
&\mathbb{E} \Bigg[ \bigg\| \int_{T_n}^t (\sigma(h,r^1)_s - \sigma(h,r^2)_s)
dW_s \bigg\|^p \Bigg] 
\\ &\leq C_p \mathbb{E} \left[ \bigg( \int_{T_n}^t \| \sigma(h,r^1)_s
- \sigma(h,r^2)_s \|_{L_2^0}^2 ds \bigg)^{\frac{p}{2}} \right]
\\ &\leq C_p (t-T_n)^{\frac{p}{2} - 1} \int_{T_n}^t
\mathbb{E} [\| \sigma(h,r^1)_s - \sigma(h,r^2)_s \|_{L_2^0}^p] ds 
\\ &\leq 2^p C_p (t-T_n)^{\frac{p}{2} - 1} \bigg( \int_{T_n}^t L(s)^p ds \bigg) \| r^1 - r^2 \|_{[T_n,T_{n+1}]}^p,
\end{align*}
with a constant $C_p > 0$. By means of the Bichteler-Jacod inequality (see \cite[Lemma 3.1]{Marinelli-Prevot-Roeckner}) and (\ref{Lip-gamma-Lp}) we get
\begin{align*}
&\mathbb{E} \Bigg[ \bigg\| \int_{T_n}^t \int_E (\gamma(h,r^1)(s,x) -
\gamma(h,r^2)(s,x)) (\mu(ds,dx) - F(dx) ds) \bigg\|^p \Bigg]
\\ &\leq N \mathbb{E} \left[ \int_{T_n}^t \int_E \| \gamma(h,r^1)(s,x) - \gamma(h,r^2)(s,x) \|^p F (dx) ds \right]
\\ &\quad + N \mathbb{E} \left[ \int_{T_n}^t \bigg(  \int_E \| \gamma(h,r^1)(s,x) - \gamma(h,r^2)(s,x) \|^2
F (dx) \bigg)^{\frac{p}{2}} ds \right]
\\ &\leq 2^{p+1} N \bigg( \int_{T_n}^t L(s)^p ds \bigg) \| r^1 - r^2 \|_{[T_n,T_{n+1}]}^p
\end{align*}
with a constant $N = N(p,t) > 0$. Proceeding as in the proof of
Lemma \ref{lemma-contraction-seq}, we obtain, after choosing an appropriate sequence $(T_n)_{n \in \mathbb{N}}$, that the fixed point mappings $\Gamma_{h}$ for $h \in C_{\rm ad}[0,T_n]$ are contractions on $C_{\rm ad}[T_n,T_{n+1}]$.

\section{Stability of stochastic differential equations}\label{sec-SDE-stability}

We shall now deal with stability of stochastic differential
equations of the kind (\ref{equation-strong}). Again these are
standard results which we do only give for the sake of completeness.
Using the method of the moving frame, we will transfer the results to stochastic partial differential equations in Section
\ref{sec-stability-SPDE}.

As in Section \ref{sec-SDE-ex}, we assume that $\alpha : C_{\rm
ad}(\mathbb{R}_+) \rightarrow H_{\mathcal{P}}$, $\sigma : C_{\rm ad}(\mathbb{R}_+)
\rightarrow (L_2^0)_{\mathcal{P}}$ and $\gamma : C_{\rm ad}(\mathbb{R}_+)
\rightarrow H_{\mathcal{P} \otimes \mathcal{E}}$ fulfill Assumptions
\ref{ass-0-strong}, \ref{ass-1-strong}, \ref{ass-2-strong}. Furthermore, let, for
each $n \in \mathbb{N}$, $\alpha_n : C_{\rm ad}(\mathbb{R}_+) \rightarrow
H_{\mathcal{P}}$, $\sigma_n : C_{\rm ad}(\mathbb{R}_+) \rightarrow
(L_2^0)_{\mathcal{P}}$ and $\gamma_n : C_{\rm ad}(\mathbb{R}_+) \rightarrow
H_{\mathcal{P} \otimes \mathcal{E}}$ be given. We make the following
additional assumptions.

\begin{assumption}\label{ass-3a-strong}
We assume that for all $T \in \mathbb{R}_+$ and all $r^1,r^2 \in C_{\rm ad}(\mathbb{R}_+)$ with $r^1|_{[0,T]} = r^2|_{[0,T]}$ we have
\begin{align*}
\alpha_n(r^1)|_{[0,T]} &= \alpha_n(r^2)|_{[0,T]}, \quad n \in \mathbb{N}
\\ \sigma_n(r^1)|_{[0,T]} &= \sigma_n(r^2)|_{[0,T]}, \quad n \in \mathbb{N}
\\ \gamma_n(r^1)|_{[0,T] \times E} &= \gamma_n(r^2)|_{[0,T] \times E}, \quad n \in \mathbb{N}.
\end{align*}
\end{assumption}

\begin{assumption}\label{ass-3-strong}
Denoting by $\mathbf{0} \in C_{\rm ad}(\mathbb{R}_+)$ the zero process, we
assume that
\begin{align*}
t \mapsto \mathbb{E} [ \| \alpha_n({\mathbf{0}})_t \|^2 ] &\in
\mathcal{L}_{\rm loc}^{1}(\mathbb{R}_+), \quad n \in \mathbb{N}
\\ t \mapsto \mathbb{E} [ \| \sigma_n({\mathbf{0}})_t \|_{L_2^0}^2 ] &\in \mathcal{L}_{\rm loc}^{1}(\mathbb{R}_+), \quad n \in \mathbb{N}
\\ t \mapsto \mathbb{E} \bigg[ \int_E \| \gamma_n({\mathbf{0}})(t,x)
\|^2 F(dx) \bigg] &\in \mathcal{L}_{\rm loc}^{1}(\mathbb{R}_+), \quad n \in \mathbb{N}.
\end{align*}
\end{assumption}

\begin{assumption}\label{ass-4-strong}
We assume that for all $t \in \mathbb{R}_+$ we have
\begin{align}\label{Lip-alpha-SDE-strong-n}
\mathbb{E} [ \| \alpha_n(r^1)_t - \alpha_n(r^2)_t \|^2 ] &\leq L(t)^2 \| r^1 - r^2
\|_t^2,
\\ \label{Lip-sigma-SDE-strong-n} \mathbb{E} [ \|
\sigma_n(r^1)_t - \sigma_n(r^2)_t \|_{L_2^0}^2 ] &\leq L(t)^2 \| r^1 - r^2
\|_t^2,
\\ \label{Lip-gamma-SDE-strong-n} \mathbb{E} \bigg[ \int_E \| \gamma_n(r^1)(t,x)
- \gamma_n(r^2)(t,x) \|^2 F(dx) \bigg] &\leq L(t)^2 \| r^1 - r^2
\|_t^2
\end{align}
for all $r^1,r^2 \in C_{\rm ad}(\mathbb{R}_+)$ and $n \in \mathbb{N}$, where $L
\in \mathcal{L}_{\rm loc}^{2}(\mathbb{R}_+)$ denotes the function from
Assumption \ref{ass-2-strong}
\end{assumption}

\begin{remark}
Notice the slight difference of the previous Assumption \ref{ass-4-strong} to Assumption \ref{ass-2-strong} for each $ \alpha_n $, $ \sigma_n $ and $ \gamma_n $, namely, that the function $ L $ does not depend on $ n \in \mathbb{N} $.
\end{remark}

Furthermore, let $t_0 \in \mathbb{R}_+$, $h \in \mathcal{C}_{\rm ad}[0,t_0]$
and for each $n \in \mathbb{N}$ let $h^n \in
\mathcal{C}_{\rm ad}[0,t_0]$ and $B_n \in \mathcal{E}$ be
given.

According to Theorem \ref{thm-existence}, there exists a unique
solution $r \in \mathcal{C}_{\rm ad}(\mathbb{R}_+)$ for
(\ref{equation-strong}) with $r|_{[0,t_0]} = h$ with c\`adl\`ag paths on $[t_0,\infty)$ satisfying
(\ref{solution-in-S2}), and for each $n \in \mathbb{N}$ there exists
a unique solution $r^n \in \mathcal{C}_{\rm ad}(\mathbb{R}_+)$ for
\begin{align*}
\left\{
\begin{array}{rcl}
dr_t^n & = & \alpha_n(r^n)_t dt + \sigma_n(r^n)_t dW_t + \int_{B_n}
\gamma_n(r^n)(t,x) (\mu(dt,dx) - F(dx) dt) \medskip
\\ r^n|_{[0,t_0]} & = & h^n,
\end{array}
\right.
\end{align*}
with c\`adl\`ag paths on $[t_0,\infty)$
satisfying $\mathbb{E} [ \sup_{t \in [t_0,T]} \| r_t^n \|^2 ] <
\infty$ for all $T \geq t_0$. 

We also make the following assumption, in which $r \in \mathcal{C}_{\rm
ad}(\mathbb{R}_+)$ denotes the solution for (\ref{equation-strong}) with $r|_{[0,t_0]} =
h$.

\begin{assumption}\label{ass-convergence}
We assume that $B_n \uparrow E$ and
\begin{align*}
\alpha_n([r]) \mathbbm{1}_{[t_0,\infty)} &\rightarrow \alpha([r]) \mathbbm{1}_{[t_0,\infty)} \quad \text{in $L^2(\lambda;H)$,}
\\ \sigma_n([r]) \mathbbm{1}_{[t_0,\infty)} &\rightarrow \sigma([r]) \mathbbm{1}_{[t_0,\infty)} \quad
\text{in $L^2(W;L_2^0)$,}
\\ \gamma_n([r]) \mathbbm{1}_{[t_0,\infty)} &\rightarrow \gamma([r]) \mathbbm{1}_{[t_0,\infty)} \quad \text{in
$L^2(\mu;H)$.}
\end{align*}
\end{assumption}

Notice that, by Assumption \ref{ass-convergence}, for all $T \geq t_0$ we have
\begin{equation}\label{def-Cn}
\begin{aligned}
C_n(T,r) &:= \bigg( \mathbb{E}\bigg[ \int_{t_0}^T \| \alpha([r])_s -
\alpha_n([r])_s \|^2 ds \bigg] + \mathbb{E}\bigg[
\int_{t_0}^T \| \sigma([r])_s - \sigma_n([r])_s \|_{L_2^0}^2 ds \bigg]
\\ &\quad + \mathbb{E} \bigg[ \int_{t_0}^T \int_E \| \gamma([r])(s,x) -
\gamma_n([r])(s,x) \|^2 F(dx) ds \bigg] 
\\ &\quad + \mathbb{E}\bigg[ \int_{t_0}^T \int_{E \setminus B_n} \| \gamma([r])(s,x)
\|^2 F(dx)ds \bigg] \bigg)^{\frac{1}{2}} \rightarrow 0 \quad \text{as $n \rightarrow \infty$.}
\end{aligned}
\end{equation}
For a compact interval $I \subset \mathbb{R}_+$
we shall also consider the norm
\begin{align*}
\| r \|_{S^2(I)} := \sqrt{\mathbb{E} \bigg[ \sup_{t \in I} \| r_t \|^2 \bigg]}.
\end{align*}
By Theorem \ref{thm-existence}, for any $T \geq t_0$ we have $\| r \|_{S^2[t_0,T]} < \infty$, where $r \in \mathcal{C}_{\rm ad}(\mathbb{R}_+)$ denotes the solution for (\ref{equation-strong}) with $r|_{[0,t_0]} = h$.

\begin{proposition}\label{prop-stability-strong}
Suppose that Assumptions \ref{ass-0-strong}, \ref{ass-1-strong}, \ref{ass-2-strong}, \ref{ass-3a-strong}, \ref{ass-3-strong}, \ref{ass-4-strong} and \ref{ass-convergence} are
fulfilled. Then, there exist maps $K_1,K_2 : \mathbb{R}_+ \rightarrow \mathbb{R}_+$, only depending on the Lipschitz function $L$, such that the following statements are valid:

\begin{enumerate}
\item If $h^n \rightarrow h$ in $C_{\rm ad}[0,t_0]$, then for each $T \geq t_0$ we have the estimate
\begin{align}\label{est-expectation}
\sup_{t \in [0,T]} \mathbb{E}[ \| r_t - r_t^n \|^2 ] &\leq K_1 \big( \| h - h^n \|_{[0,t_0]}^2 + C_n^2 \big)
\rightarrow 0 \quad \text{for $n \rightarrow \infty$,}
\end{align}
where $K_1 = K_1(T)$ and $C_n = C_n(T,r)$ is defined in (\ref{def-Cn}).

\item If even $h^n \rightarrow h$ in $S^2[0,t_0]$, then for each $T \geq t_0$ we have the estimate
\begin{align}\label{est-sup-expectation}
\mathbb{E}\bigg[ \sup_{t \in [0,T]} \| r_t - r_t^n \|^2 \bigg] &\leq
K_2 \big( \| h - h^n \|_{S^2[0,t_0]}^2 + C_n^2 \big)
\rightarrow 0 \quad \text{for $n \rightarrow \infty$,}
\end{align}
where $K_2 = K_2(T)$ and $C_n = C_n(T,r)$ is defined in (\ref{def-Cn}).
\end{enumerate}
\end{proposition}

\begin{proof}
Let $T \geq t_0$ and $n \in \mathbb{N}$ be arbitrary. By
H\"older's inequality, the It\^o-isometries (\ref{Ito-isometry-Wiener}),
(\ref{Ito-isometry-mu}) and the Lipschitz conditions
(\ref{Lip-alpha-SDE-strong-n}), (\ref{Lip-sigma-SDE-strong-n}),
(\ref{Lip-gamma-SDE-strong-n}) we obtain, by writing
\begin{align*}
&\int_{t_0}^t \int_E \gamma([r])(s,x)(\mu(ds,dx) - F(dx)ds) 
\\ &\quad - \int_{t_0}^t \int_{B_n} \gamma_n([r^n])(s,x)(\mu(ds,dx) - F(dx)ds)
\\ &= \int_{t_0}^t \int_{B_n} (\gamma([r])(s,x) - \gamma_n([r])(s,x))(\mu(ds,dx) - F(dx)ds)
\\ &\quad + \int_{t_0}^t \int_{E \setminus B_n} \gamma([r])(s,x)(\mu(ds,dx) - F(dx)ds)
\\ &\quad + \int_{t_0}^t \int_{B_n} (\gamma_n([r])(s,x) - \gamma_n([r^n])(s,x)(\mu(ds,dx) - F(dx)ds),
\end{align*}
for all $t \in [t_0,T]$ the estimate
\begin{align*}
&\| r - r^n \|_{[0,t]}^2 = \sup_{s \in [0,t]} \mathbb{E}[\| r_s - r_s^n \|^2] \leq 8 \big( \| h - h^n \|_{[0,t_0]}^2 + ((t-t_0) \vee 1) C_n(t,r)^2 \big)
\\ &\quad + 8 \sup_{s \in
[t_0,t]} \mathbb{E}\left[ \bigg\| \int_{t_0}^s (\alpha_n([r])_v - \alpha_n([r^n])_v)dv
\bigg\|^2 \right]
\\ &\quad + 8 \sup_{s \in
[t_0,t]} \mathbb{E}\left[ \bigg\| \int_{t_0}^s (\sigma_n([r])_v - \sigma_n([r^n])_v)dW_v
\bigg\|^2 \right]
\\ &\quad + 8 \sup_{s \in [t_0,t]} \mathbb{E}\left[ \bigg\| \int_{t_0}^s \int_{B_n} (\gamma_n([r])(v,x) -
\gamma_n([r^n])(v,x))(\mu(dv,dx) - F(dx) dv) \bigg\|^2
\right]
\\ &\leq 8 \big( \| h - h^n \|_{[0,t_0]}^2 + ((T-t_0) \vee 1) C_n(T,r)^2 \big)
\\ &\quad + 8 (T-t_0+2) \int_{t_0}^t L(s)^2 \| r - r^n \|_{[0,s]}^2 ds.
\end{align*}
Applying the Gronwall Lemma gives us
\begin{align*}
&\sup_{s \in [0,t]} \mathbb{E}[\| r_s - r_s^n \|^2] = \| r - r^n \|_{[0,t]}^2 
\\ &\leq 8 \big( \| h - h^n \|_{[0,t_0]}^2 + ((T-t_0) \vee 1) C_n(T,r)^2 \big) e^{8(T-t_0+2) \int_{t_0}^t L(s)^2 ds}
\end{align*}
for all $t \in [t_0,T]$, implying (\ref{est-expectation}). Analogously, by also taking into account Doob's martingale inequality \cite[Thm.
3.8]{Da_Prato}, we obtain
\begin{align*}
&\mathbb{E} \bigg[ \sup_{t \in
[0,T]} \| r_t - r_t^n \|^2 \bigg] \leq 8 \big( \| h-h^n \|_{S^2[0,t_0]}^2 + ((T-t_0) \vee 4)C_n(T,r)^2 \big)
\\ &\quad + 8 \mathbb{E}\left[ \sup_{t \in
[t_0,T]} \bigg\| \int_{t_0}^t (\alpha_n([r])_s - \alpha_n([r^n])_s)ds
\bigg\|^2 \right]
\\ &\quad + 8 \mathbb{E}\left[ \sup_{t \in
[t_0,T]} \bigg\| \int_{t_0}^t (\sigma_n([r])_s - \sigma_n([r^n])_s)dW_s
\bigg\|^2 \right]
\\ &\quad + 8 \mathbb{E}\left[ \sup_{t \in [t_0,T]} \bigg\| \int_{t_0}^t \int_{B_n} (\gamma_n([r])(s,x) -
\gamma_n([r^n])(s,x))(\mu(ds,dx) - F(dx) ds) \bigg\|^2
\right]
\\ &\leq 8 \big( \| h-h^n \|_{S^2[0,t_0]}^2 + ((T-t_0) \vee 4)C_n(T,r)^2 \big) 
\\ &\quad + 8 (T-t_0+8) \bigg( \int_{t_0}^T L(s)^2 ds \bigg) \| r - r^n \|_{[0,T]}^2.
\end{align*}
Noting that $\| h-h^n \|_{[0,t_0]} \leq \| h-h^n \|_{S^2[0,t_0]}$, inserting (\ref{est-expectation}) shows (\ref{est-sup-expectation}).
\end{proof}

\begin{remark}\label{remark-solution-map}
Fix a finite time $T \geq t_0$ and denote for $h \in \mathcal{C}_{\rm ad}[0,t_0]$ by $r^h$ the unique solution for (\ref{equation-strong}) with $r|_{[0,t_0]} = h$, which has c\`{a}dl\`{a}g paths on $[t_0,\infty)$. Restricting it to the interval $[0,T]$, estimates (\ref{est-expectation}), (\ref{est-sup-expectation}) show that the solution map $h \mapsto r^h$ is Lipschitz continuous with a constant $L = L(T) > 0$, if considered as a map $C_{\rm ad}[0,t_0] \rightarrow C_{\rm ad}[0,T]$ or as a map $S^2[0,t_0] \rightarrow S^2[0,T]$. In particular, there exists a constant $C = C(T) > 0$ such that
\begin{align*}
\sup_{t \in [0,T]} \mathbb{E} [\| r_t^h \|^2] &\leq C \Big( 1 + \sup_{t \in [0,t_0]} \mathbb{E} [ \| h_t \|^2 ] \Big), \quad h \in \mathcal{C}_{\rm ad}[0,t_0]
\\ \mathbb{E} \bigg[ \sup_{t \in [0,T]} \| r_t^h \|^2 \bigg] &\leq C \bigg( 1 + \mathbb{E}\bigg[ \sup_{t \in [0,t_0]} \| h_t \|^2 \bigg] \bigg), \quad h \in \mathcal{S}^2[0,t_0].
\end{align*}
Notice further for $t_0 = 0$ the coincidence $\mathcal{C}_{\rm ad}[0,t_0] = \mathcal{S}^2[0,t_0] = \mathcal{L}^2(\Omega,\mathcal{F}_0,\mathbb{P};H)$.
\end{remark}

\begin{remark}
Using Burkholder-Davis-Gundy and Bichteler-Jacod type arguments as in the previous section, we can, in an analogous fashion, derive the $L^p$-version of the stability result above.
\end{remark}

\section{Regular dependence on initial data for stochastic differential equations}\label{sec-SDE-regularity}

In this section, we study regular dependence on initial data for
SDEs. Some related ideas can be found in \cite{Mandrekar-Ruediger}.
By the method of the moving frame, which we present in Section
\ref{sec-SPDE-ex}, we can transfer the upcoming results to SPDEs.

We understand the question of regular dependence on initial data as a conclusion of the stability results of Section \ref{sec-SDE-stability}. We consider
\begin{align}\label{equation-strong-regularity}
\left\{
\begin{array}{rcl}
dr_t & = & \alpha(r)_t dt + \sigma(r)_t dW_t + \int_E \gamma(r)(t,x)
(\mu(dt,dx) - F(dx) dt) \medskip
\\ r|_{[0,t_0]} & = & h,
\end{array}
\right.
\end{align}
under Assumptions \ref{ass-0-strong}, \ref{ass-1-strong}, \ref{ass-2-strong}, such that we can conclude the existence and uniqueness of solutions for $t_0 \in \mathbb{R}_+$ and $ h \in C_{\rm ad}[0,t_0] $. Motivated by ideas from convenient analysis, see \cite{KriMic:97}, we fix a curve of initial data $ \epsilon \mapsto c(\epsilon) \in C_{\rm ad}[0,t_0] $, which is differentiable for all $ \epsilon $ with derivative $ c'(\epsilon) \in C_{\rm ad}[0,t_0] $. We consider the following system of equations,
\begin{align}
\left\{
\begin{array}{rcl}
dr^{\epsilon}_t & = &\alpha(r^{\epsilon})_t dt + \sigma(r^{\epsilon})_t dW_t + \int_E \gamma(r^{\epsilon})(t,x)
(\mu(dt,dx) - F(dx) dt),
\\ r^{\epsilon}|_{[0,t_0]} & =&  c(\epsilon),
\\ d\frac{r^{\epsilon}_t - r^0_t}{\epsilon}  & = & \frac{\alpha(r^{\epsilon})_t - \alpha(r^{0})_t}{\epsilon}  dt + \frac{\sigma(r^{\epsilon})_t - \sigma(r^{0})_t}{\epsilon}  dW_t + \\& & +\int_E \frac{\gamma(r^{\epsilon})(t,x) - \gamma(r^{0})(t,x)}{\epsilon} (\mu(dt,dx) - F(dx) dt),
\\ \frac{r^{\epsilon} - r^0}{\epsilon}|_{[0,t_0]} & = & \frac{c(\epsilon) - c(0)}{\epsilon},
\end{array}
\right.
\end{align}
for $ \epsilon \neq 0 $, where $r^0$ denotes the solution for (\ref{equation-strong-regularity}) with $h = c(0)$. We can consider those equations indeed as two SDEs in our sense. More precisely let
\begin{align}\label{two-SDEs-reg}
\left\{
\begin{array}{rcl}
dr^{\epsilon}_t & = &\alpha(r^{\epsilon})_t dt + \sigma(r^{\epsilon})_t dW_t + \int_E \gamma(r^{\epsilon})(t,x)
(\mu(dt,dx) - F(dx) dt),
\\ r^{\epsilon}|_{[0,t_0]} & =&  c(\epsilon),
\\ d \Delta^{\epsilon}_t  & = & \frac{\alpha(\epsilon \Delta^{\epsilon} + r^0)_t - \alpha(r^{0})_t}{\epsilon}  dt + \frac{\sigma(\epsilon \Delta^{\epsilon} + r^0)_t - \sigma(r^{0})_t}{\epsilon}  dW_t + \\& & +\int_E \frac{\gamma(\epsilon \Delta^{\epsilon} + r^0)(t,x) - \gamma(r^{0})(t,x)}{\epsilon} (\mu(dt,dx) - F(dx) dt),
\\ \Delta^{\epsilon}|_{[0,t_0]} & = & \frac{c(\epsilon) - c(0)}{\epsilon},
\end{array}
\right.
\end{align}
for $ \epsilon \neq 0 $, then this system of equations can be seen as two stochastic differential equations. We can readily check that the Assumptions \ref{ass-0-strong}, \ref{ass-1-strong}, \ref{ass-2-strong} are true for the second SDE in (\ref{two-SDEs-reg}) for every $ \epsilon \neq 0 $. Its solution is given by
\begin{align}
\Delta^{\epsilon}_t = \frac{r^{\epsilon}_t - r^0_t}{\epsilon}, \quad t \geq 0.
\end{align}
We assume now that the maps $ \alpha $, $ \sigma $ and $ \gamma $
admit directional derivatives in all directions of $ C_{\rm ad}(\mathbb{R}_+)
$. We denote those directional
derivatives at the point $ r \in C_{\rm ad}(\mathbb{R}_+) $ into direction $ v
\in C_{\rm ad}(\mathbb{R}_+) $  by $ D \alpha (r) \bullet v  $, $ D \sigma(r)
\bullet v $ and $ D \gamma(r) \bullet v $. By $D$, we always
mean the Fr\'echet derivative.

\begin{assumption}\label{ass-convergence-stability}
We define the first variation process $  J(r) \bullet w $ in
direction $ w $, where $ w \in C_{\rm ad}[0,t_0] $, to be the unique solution
of the SDE
\begin{align}\label{first-variation-equation}
\left\{
\begin{array}{rcl}
d (J(r) \bullet w)_t  & = & {\bigl(D \alpha(r) \bullet (J(r) \bullet w) \bigr)}_t  dt + {\bigl(D \sigma(r) \bullet (J(r) \bullet w) \bigr)}_t  dW_t + \\& & +\int_E {\bigl(D \gamma(r) \bullet (J(r) \bullet w) \bigr)}(t,x) (\mu(dt,dx) - F(dx) dt),
\\ (J(r) \bullet w)|_{[0,t_0]} & = & w,
\end{array}
\right.
\end{align}
where $ r $ solves equation \eqref{equation-strong-regularity}. We assume that Assumptions \ref{ass-0-strong}, \ref{ass-1-strong}, \ref{ass-2-strong} are true for equation \eqref{first-variation-equation}. We assume furthermore that
\begin{align*}
\frac{\alpha(\epsilon (J(r)\bullet w) + r) - \alpha(r)}{\epsilon} \mathbbm{1}_{[t_0,\infty)} & \rightarrow
D \alpha(r) \bullet (J(r) \bullet w) \mathbbm{1}_{[t_0,\infty)},  \\
\frac{\sigma(\epsilon (J(r)\bullet w) + r) - \sigma(r)}{\epsilon} \mathbbm{1}_{[t_0,\infty)} & \rightarrow D \sigma(r) \bullet (J(r) \bullet w) \mathbbm{1}_{[t_0,\infty)}, \\ \frac{\gamma(\epsilon (J(r) \bullet w)+ r) - \gamma(r)}{\epsilon} \mathbbm{1}_{[t_0,\infty)} & \rightarrow D \gamma(r) \bullet (J(r) \bullet w) \mathbbm{1}_{[t_0,\infty)}
\end{align*}
as $ \epsilon \to 0 $ in the respective spaces $L^2(\lambda;H)$, $L^2(W;L_2^0)$ and $L^2(\mu;H)$. The process $ r $
denotes the solution of equation \eqref{equation-strong-regularity}
and $ J(r) \bullet w $ denotes the solution of the first variation
equation \eqref{first-variation-equation}.
\end{assumption}

\begin{proposition}\label{prop-regularity-strong}
Suppose that Assumptions \ref{ass-3a-strong}, \ref{ass-3-strong}, \ref{ass-4-strong} for equation
\begin{align}
\left\{
\begin{array}{rcl}
d \Delta^{\epsilon}_t  & = & \frac{\alpha(\epsilon \Delta^{\epsilon} + r^0)_t - \alpha(r^{0})_t}{\epsilon}  dt + \frac{\sigma(\epsilon \Delta^{\epsilon} + r^0)_t - \sigma(r^{0})_t}{\epsilon}  dW_t + \\& & +\int_E \frac{\gamma(\epsilon \Delta^{\epsilon} + r^0)(t,x) - \gamma(r^{0})(t,x)}{\epsilon} (\mu(dt,dx) - F(dx) dt),
\\ \Delta^{\epsilon}|_{[0,t_0]} & = & \frac{c(\epsilon) - c(0)}{\epsilon},
\end{array}
\right.
\end{align}
are valid in the obvious sense for $ \epsilon \neq 0 $ in a neighborhood of $ 0 $, and assume that Assumption \ref{ass-convergence-stability} is fulfilled for $ w = c'(0) $ for a chosen curve of initial values $ \epsilon \mapsto c(\epsilon) $. Then, for each $T \geq t_0$ we have the estimate
\begin{align}\label{est-expectation-regularity}
\sup_{t \in [0,T]} \mathbb{E}[ \| {(J(r) \bullet w)}_t - \Delta^{\epsilon}_t \|^2 ] \leq K_1 \Bigg( \bigg\| c'(0) - \frac{c(\epsilon) - c(0)}{\epsilon} \bigg\|_{[0,t_0]}^2 + C_{\epsilon}^2 \Bigg)
\rightarrow 0
\end{align}
for $\epsilon \rightarrow 0$,
and if $\epsilon \mapsto c(\epsilon)$ is even a curve in $S^2[0,t_0]$, then for each $T \geq t_0$ we have the estimate
\begin{align}\label{est-sup-expectation-regularity}
\mathbb{E}\bigg[ \sup_{t \in [0,T]} \|  {(J(r) \bullet w)}_t - \Delta^{\epsilon}_t \|^2 \bigg] \leq
K_2 \Bigg( \bigg\| c'(0) - \frac{c(\epsilon) - c(0)}{\epsilon} \bigg\|_{S^2[0,t_0]}^2 + C_{\epsilon}^2 \Bigg)
\rightarrow 0
\end{align}
for $\epsilon \rightarrow 0$.
In particular, the map $ w \mapsto J(r) \bullet w $ is linear and continuously depending on $ w $ in the sense that
for every $T \geq t_0$ we have
\begin{align}\label{est-expectation-firstvar}
\sup_{t \in [0,T]} \mathbb{E} [ \|  {(J(r) \bullet w)}_t -
{(J(r) \bullet w^n)}_t \|^2 ] \leq K_1 \| w - w^n \|_{[0,t_0]}^2
\rightarrow 0
\end{align}
for variation of the initial value $ w^n \rightarrow w \in C_{\rm ad}[0,t_0] $, and if $\epsilon \mapsto c(\epsilon)$ is even a curve in $S^2[0,t_0]$, then for every $T \geq t_0$ we have
\begin{align}\label{est-sup-expectation-firstvar}
\mathbb{E}\bigg[ \sup_{t \in [0,T]} \|  {(J(r) \bullet w)}_t -
{(J(r) \bullet w^n)}_t \|^2 \bigg] \leq K_2 \| w - w^n \|_{S^2[0,t_0]}^2
\rightarrow 0
\end{align}
for variation of the initial value $ w^n \rightarrow w \in S^2[0,t_0] $.
\end{proposition}

\begin{remark}
The notion $ C_{\epsilon} $ is defined corresponding to (\ref{def-Cn}) and $K_1,K_2$ according to Proposition \ref{prop-stability-strong}.
\end{remark}

\begin{proof}
The assertion is a corollary of Proposition \ref{prop-stability-strong}. Assumption \ref{ass-convergence-stability} corresponds precisely to Assumption \ref{ass-convergence}, which is needed for the proof of Proposition \ref{prop-stability-strong}. Note that the ``continuous'' parameter $ \epsilon $ replaces the index $ n $, which does not cause any problems, since we do not speak about almost sure convergence results here.
\end{proof}

\begin{remark}
Fix a finite time $T \geq t_0$ and a curve of initial data $\epsilon \mapsto c(\epsilon) \in C_{\rm ad}[0,t_0]$ or $\epsilon \mapsto c(\epsilon) \in S^2[0,t_0]$. Then, we can consider the curve of solution processes $\epsilon \mapsto r^{\epsilon} \in C_{\rm ad}[0,T]$ or $\epsilon \mapsto r^{\epsilon} \in S^2[0,T]$, respectively. By Remark \ref{remark-solution-map} we already know that the solution map $\epsilon \mapsto r^{\epsilon}$ is continuous. Now, estimates (\ref{est-expectation-regularity}), (\ref{est-sup-expectation-regularity}) show that, subject to our previous assumptions, $\epsilon \mapsto r^{\epsilon}$ is also differentiable with derivative $\epsilon \mapsto J(r) \bullet c'(\epsilon)$. Moreover, regarding the variation $w \mapsto J(r) \bullet w$ of the initial value as a linear map $C_{\rm ad}[0,t_0] \rightarrow C_{\rm ad}[0,T]$ or as a linear map $S^2[0,t_0] \rightarrow S^2[0,T]$, estimates (\ref{est-expectation-firstvar}), (\ref{est-sup-expectation-firstvar}) show its continuity.
\end{remark}

Considering the construction for all possible curves of initial values $ c $ we can define the first (and possibly higher) variation processes in a coherent way for all variations of the initial values and also for variations of the process up to time $ t $ by shifting $ \mathcal{F}_t $ to $ \mathcal{F}_0 $. Properties of this variation process can be established by considering the equation, which follows right from Proposition \ref{prop-regularity-strong},
\begin{align}
r^{\epsilon} - r^0 = \int_0^{\epsilon} J(r) \bullet c'(\eta) \, d \eta
\end{align}
and which reveals the true meaning of the first variation process.

\section{Solution concepts for stochastic partial differential equations}\label{sec-SPDE-def}

When dealing with SPDEs there are several solution concepts, which
we will discuss in this section. The main difficulty is that
solutions of SPDEs usually leave the realm of semi-martingales and
one therefore has to modify the usual semi-martingale decomposition.
The method of the moving frame, which will be presented in the next
section, is a new approach how to handle this problem.

In this section, we review the well-known concepts of strong, weak
and mild solutions and show, how they are related. The proofs from
\cite{Da_Prato} (or \cite{P-Z-book}) can be transferred to the
present situation, whence we keep this section rather short. The
decisive tool in order to prove Lemma \ref{lemma-mild-weak} is an
appropriate Stochastic Fubini Theorem with respect to Poisson
measures, which we provide in Appendix \ref{a}.

Now let $(S_t)_{t \geq 0}$ be a $C_0$-semigroup on the separable
Hilbert space $H$ with infinitesimal generator $A : \mathcal{D}(A)
\subset H \rightarrow H$. We denote by $A^* : \mathcal{D}(A^*)
\subset H \rightarrow H$ the adjoint operator of $A$. Recall that
the domains $\mathcal{D}(A)$ and $\mathcal{D}(A^*)$ are dense in
$H$, see, e.g., \cite[Satz VII.4.6, p. 351]{Werner}.

In this section, we are interested in stochastic partial
differential equations of the form
\begin{align}\label{general-SPDE}
\left\{
\begin{array}{rcl}
dr_t & = & (Ar_t + \alpha(r)_t)dt + \sigma(r)_t dW_t + \int_E
\gamma(r)(t,x) (\mu(dt,dx) - F(dx) dt)
\medskip
\\ r|_{[0,t_0]} & = & h
\end{array}
\right.
\end{align}
where $\alpha : C_{\rm ad}(\mathbb{R}_+) \rightarrow H_{\mathcal{P}}$, $\sigma
: C_{\rm ad}(\mathbb{R}_+) \rightarrow (L_2^0)_{\mathcal{P}}$ and $\gamma :
C_{\rm ad}(\mathbb{R}_+) \rightarrow H_{\mathcal{P} \otimes \mathcal{E}}$. Fix $t_0 \in \mathbb{R}_+$ and $h \in \mathcal{C}_{\rm ad}[0,t_0]$.

\begin{definition}
A process $r \in \mathcal{C}_{\rm ad}(\mathbb{R}_+)$ is called a {\rm strong solution}
for (\ref{general-SPDE}) if we have $r|_{[0,t_0]} = h$, $\mathbb{P}(r_t \in
\mathcal{D}(A)) = 1$, $t \geq t_0$, the relations $A(r \mathbbm{1}_{[t_0,\infty)}) + \alpha([r]) \mathbbm{1}_{[t_0,\infty)} \in
\mathcal{L}_{\rm loc}^1(\lambda;H)$, $\sigma([r]) \mathbbm{1}_{[t_0,\infty)} \in \mathcal{L}_{\rm loc}^2(W;L_0^2)$,
$\gamma([r]) \mathbbm{1}_{[t_0,\infty)} \in \mathcal{L}_{\rm loc}^2(\mu;H)$ and almost surely
\begin{equation}\label{strong-solution}
\begin{aligned}
r_t &= h_{t_0} + \int_{t_0}^t (A r_s + \alpha([r])_s)ds + \int_{t_0}^t
\sigma([r])_s dW_s
\\ &\quad + \int_{t_0}^t \int_E \gamma([r])(s,x) (\mu(ds,dx) -
F(dx) ds), \quad t \geq t_0.
\end{aligned}
\end{equation}
\end{definition}

\begin{definition}
A process $r \in \mathcal{C}_{\rm ad}(\mathbb{R}_+)$ is called a {\rm weak solution} for
(\ref{general-SPDE}) if $r|_{[0,t_0]} = h$, $\alpha([r]) \mathbbm{1}_{[t_0,\infty)} \in \mathcal{L}_{\rm
loc}^1(\lambda;H)$, $\sigma([r]) \mathbbm{1}_{[t_0,\infty)} \in \mathcal{L}_{\rm loc}^2(W;L_0^2)$,
$\gamma([r]) \mathbbm{1}_{[t_0,\infty)} \in \mathcal{L}_{\rm loc}^2(\mu;H)$ and for all $\zeta \in
\mathcal{D}(A^*)$ we have almost surely
\begin{equation}\label{weak-solution}
\begin{aligned}
\langle \zeta,r_t \rangle &= \langle \zeta,h_{t_0} \rangle + \int_{t_0}^t (
\langle A^* \zeta, r_s \rangle + \langle \zeta, \alpha([r])_s \rangle)ds +
\int_{t_0}^t \langle \zeta, \sigma([r])_s \rangle dW_s
\\ &\quad + \int_{t_0}^t \int_E \langle \zeta, \gamma([r])(s,x) \rangle
(\mu(ds,dx) - F(dx) ds), \quad t \geq t_0.
\end{aligned}
\end{equation}
\end{definition}

\begin{definition}
A process $r \in \mathcal{C}_{\rm ad}(\mathbb{R}_+)$ is called a {\rm mild solution} for
(\ref{general-SPDE}) if $r|_{[0,t_0]} = h$, $\alpha([r]) \mathbbm{1}_{[t_0,\infty)} \in \mathcal{L}_{\rm
loc}^1(\lambda;H)$, $\sigma([r]) \mathbbm{1}_{[t_0,\infty)} \in \mathcal{L}_{\rm loc}^2(W;L_0^2)$,
$\gamma([r]) \mathbbm{1}_{[t_0,\infty)} \in \mathcal{L}_{\rm loc}^2(\mu;H)$ and we have almost surely
\begin{equation}\label{mild-solution}
\begin{aligned}
r_t &= S_{t-t_0} h_{t_0} + \int_{t_0}^t S_{t-s} \alpha([r])_s ds + \int_{t_0}^t
S_{t-s} \sigma([r])_s dW_s
\\ &\quad + \int_{t_0}^t \int_E S_{t-s} \gamma([r])(s,x) (\mu(ds,dx) -
F(dx) ds), \quad t \geq t_0.
\end{aligned}
\end{equation}
\end{definition}

\begin{remark}
For all the three just defined solution concepts {\em uniqueness} of solutions for (\ref{general-SPDE}) is, as in Definition \ref{def-solution-SDE},
meant up to indistinguishability on the interval $[t_0,\infty)$.
\end{remark}

\begin{lemma}
Let $r \in \mathcal{C}_{\rm ad}(\mathbb{R}_+)$ be a strong solution for
(\ref{general-SPDE}). Then, $r$ is also a weak
solution for (\ref{general-SPDE}).
\end{lemma}

\begin{proof}
For all $\zeta \in \mathcal{D}(A^*)$ we have
\begin{align*}
\langle \zeta,r_t \rangle &= \langle \zeta,h_{t_0} \rangle + \int_{t_0}^t
\langle \zeta, Ar_s + \alpha([r])_s \rangle ds + \int_{t_0}^t \langle
\zeta, \sigma([r])_s \rangle d W_s
\\ &\quad + \int_{t_0}^t \int_E \langle \zeta, \gamma([r])(s,x) \rangle
(\mu(ds,dx) - F(dx) ds), \quad t \geq t_0
\end{align*}
implying that $r$ is also a weak solution for (\ref{general-SPDE}), because $\langle \zeta, Ah \rangle = \langle A^*
\zeta, h \rangle$ for all $h \in \mathcal{D}(A)$.
\end{proof}

\begin{lemma}\label{lemma-weak-mild}
Let $r \in \mathcal{C}_{\rm ad}(\mathbb{R}_+)$ be a weak solution for
(\ref{general-SPDE}). Then, $r$ is also a mild
solution for (\ref{general-SPDE}).
\end{lemma}

\begin{proof}
Let $T \geq t_0$ be arbitrary. As in the proof of \cite[Thm.
9.15]{P-Z-book} we show that
\begin{align*}
&\langle g(t), r_t \rangle = \langle g(t_0), h_{t_0} \rangle + \int_{t_0}^t
\Big( \langle g'(s) + A^* g(s), r_s \rangle + \langle g(s),
\alpha([r])_s \rangle \Big) ds
\\ & \quad + \int_{t_0}^t \langle g(s), \sigma([r])_s \rangle d
W_s + \int_{t_0}^t \int_E \langle g(s), \gamma([r])(s,x) \rangle
(\mu(ds,dx) - F(dx) ds)
\end{align*}
for all $g \in C^1([t_0,T];\mathcal{D}(A^*))$ and $t \in [t_0,T]$. For an arbitrary $t \geq t_0$ and an arbitrary $\zeta \in \mathcal{D}(A^*)$ we apply
this identity to $g(s) := S_{t-s}^* \zeta$, $s \in [t_0,t]$, which
yields that the process $r$ is also a mild solution for (\ref{general-SPDE}).
\end{proof}

\begin{lemma}\label{lemma-mild-weak}
Let $r \in \mathcal{C}_{\rm ad}(\mathbb{R}_+)$ be a mild solution for
(\ref{general-SPDE}) such that $\sigma([r]) \mathbbm{1}_{[t_0,\infty)} \in
L^2(W;L_2^0)$ and $\gamma([r]) \mathbbm{1}_{[t_0,\infty)} \in L^2(\mu;H)$. Then, $r$ is also a
weak solution for (\ref{general-SPDE}).
\end{lemma}

\begin{proof}
We proceed as in the proof of \cite[Thm. 9.15]{P-Z-book}. The change
of order of integration for the stochastic integrals with respect to
the compensated Poisson random measure is valid by the Stochastic
Fubini Theorem \ref{thm-Fubini-Poisson} provided in Appendix
\ref{a}.
\end{proof}

\section{Existence and uniqueness of mild and weak solutions for
stochastic partial differential equations}\label{sec-SPDE-ex}

In this section we introduce the method of the moving frame, which has been announced in the introduction. 
Loosely speaking we apply a time-dependent coordinate transformation to the SPDE such that ``from the point of view of
the moving frame'' the SPDE looks like an SDE with appropriately transformed coefficients. 
The method is in contrast to the point of view, that an SPDE is a PDE together with a non-linear stochastic perturbation. 
Here we consider an SPDE rather as a time-transformed SDE, where the time transform contains the respective PDE aspect.

We apply this method for an ``easy'' proof of existence and uniqueness in this general setting. The key argument, which allows to apply the method, is the Sz\H{o}kefalvi-Nagy theorem, which has been brought to our attention by \cite{Seidler}. We emphasize that in our article we do not need a particular representation of the Hilbert space involved in the Sz\H{o}kefalvi-Nagy theorem (see the subsequent remark). The Sz\H{o}kefalvi-Nagy theorem is a ``ladder'', which allows us to ``climb'' towards several new assertions, but which is not necessary to understand the statements of those assertions.

During this section, we impose the
following assumption.

\begin{assumption}\label{ass-group}
There exist another separable Hilbert space $\mathcal{H}$, a
$C_0$-group $(U_t)_{t \in \mathbb{R}}$ on $\mathcal{H}$ and
continuous linear operators $\ell \in L(H,\mathcal{H})$, $\pi \in
L(\mathcal{H},H)$ such that the diagram
\[ \begin{CD}
\mathcal{H} @>U_t>> \mathcal{H}\\
@AA\ell A @VV\pi V\\
H @>S_t>> H
\end{CD} \]
commutes for every $t \in \mathbb{R}_+$, that is
\begin{align}\label{diagram-commutes}
\pi U_t \ell = S_t \quad \text{for all $t \in \mathbb{R}_+$.}
\end{align}
In particular, we see that $ \pi \ell = {\rm Id} $.
\end{assumption}

\begin{remark}
In the spirit of \cite{Nagy}, the group $(U_t)_{t \in \mathbb{R}}$
is a {\rm dilation} of the semigroup $(S_t)_{t \geq 0}$.
\end{remark}

\begin{remark}\label{remark-H}
Assumption \ref{ass-group} is not only frequently fulfilled, which
seems surprising at a first view, but it is also possible to
describe the respective Hilbert space $ \mathcal{H} $ more
precisely. Take for instance a self-adjoint strongly continuous
semigroup of contractions $ S $ on the complex Hilbert space $ H $,
then -- as a part of the Sz\H{o}kefalvi-Nagy theorem -- the map $ t
\mapsto S_{| t |} $, where the semigroup is extended by $S_{-t} :=
S_t$ for $t \geq 0$, is a strongly continuous, positive definite
map, i.e.~for all $ \psi_1,\ldots,\psi_n \in H $ and all real times
$ t_1,\ldots,t_n $ the matrix $ \bigl ( \langle S_{| t_i - t_j |}
\psi_i, \psi_j \rangle \bigr ) $ is positive definite. A positive
definite map with values in bounded linear operators can be
considered as characteristic function of a vector-valued measure $
\eta $ taking values in positive operators on $ H $. One can define
the Hilbert space $ \mathcal{H} = L^2(\mathbb{R},\eta;H) $, i.e.~the
space of square-integrable $H$-valued measurable maps $ f $, such
that the integral
$$
\int_{\mathbb{R}} \langle f(x) , \eta(dx) f(x) \rangle < \infty
$$
is finite. $ H $ can be embedded via the constant maps $ f(x) \equiv h $ for $ h \in H $ and $ x \in \mathbb{R} $ and the semigroup $ U $ is defined via
$$
U_t f (x ) = \exp(itx) f(x)
$$
for $ t ,x \in \mathbb{R} $. Consequently, more precise analysis of the respective generator of $ S $ on $ \mathcal{H} $ can be performed. Details of the previous considerations and impacts on SPDEs will be presented elsewhere.
\end{remark}

According to Proposition \ref{prop-group} below, Assumption
\ref{ass-group} is in particular satisfied if the semigroup
$(S_t)_{t \geq 0}$ is pseudo-contractive.

\begin{definition}
The $C_0$-semigroup $(S_t)_{t \geq 0}$ is called {\rm
pseudo-contractive} if there exists $\omega \in \mathbb{R}$ such
that
\begin{align}\label{def-pseudo-contr}
\| S_t \| \leq e^{\omega t}, \quad t \geq 0.
\end{align}
\end{definition}

\begin{remark}
Sometimes in the literature, e.g., see \cite{Marinelli-Prevot-Roeckner}, the notion quasi-contractive is used instead of pseudo-contractive.
\end{remark}

\begin{remark}\label{remark-dissipative}
By the theorem of Lumer-Phillips, a densely defined operator $A$ generates a pseudo-contractive
semigroup $(S_t)_{t \geq 0}$ with growth estimate (\ref{def-pseudo-contr})
for some $\omega \geq 0$ if and only if $A$ is
$\omega$-m-dissipative, that is, $A - \omega$ is dissipative, which
means
\begin{align}\label{dissipative}
\langle Ah,h \rangle \leq \omega \| h \|^2 \quad \text{for all $h
\in \mathcal{D}(A)$,}
\end{align}
and there exists $\lambda > 0$ such that $\lambda + \omega - A$ is
surjective. For example, consider the Hilbert space $H = L^2(0,\infty)$ and the Laplace operator $A = \Delta$ defined by $\Delta h = h''$ on the
Sobolev space $\mathcal{D}(\Delta) = H_0^1(0,\infty) \cap W^2(0,\infty)$. Then, $\Delta$ is densely defined, because $C_0^{\infty}(0,\infty)$ is dense in $L^2(0,\infty)$. Let us check the dissipativity of
$\Delta$. For $h \in H_0^1(0,\infty) \cap W^2(0,\infty)$ choose a sequence $(\varphi_n)_{n \in \mathbb{N}} \subset C_0^{\infty}(0,\infty)$ with $\| h - \varphi_n \|_{H_0^1} \rightarrow 0$.
By integration by parts, we have
\begin{align*}
\langle h'',h \rangle_{L^2} = \lim_{n \rightarrow \infty} \langle h'', \varphi_n \rangle_{L^2} = \lim_{n \rightarrow \infty} \langle h', \varphi_n' \rangle_{L^2} = \langle h',h' \rangle_{L^2} \leq 0,
\end{align*}
showing (\ref{dissipative}) with $\omega = 0$. For $\lambda > 0$ and $f \in L^2(0,\infty)$ there exists a unique solution $h \in H_0^1(0,\infty) \cap W^2(0,\infty)$ of the second order differential equation
\begin{align*}
\lambda h - \Delta h = f,
\end{align*}
see \cite[Thm. 8.2.7]{Krylov}. Hence, $\lambda - \Delta$ is surjective.
\end{remark}

For every $C_0$-semigroup $(S_t)_{t \geq 0}$ there exist constants
$M \geq 1$ and $\omega \in \mathbb{R}$ such that
\begin{align}\label{semigroup-growth}
\| S_t \| \leq M e^{\omega t}, \quad t \geq 0
\end{align}
see, e.g., \cite[Lemma VII.4.2]{Werner}. Hence, in other words, the
semigroup $(S_t)_{t \geq 0}$ is contractive if we can choose $M = 1$
and $\omega = 0$ in (\ref{semigroup-growth}), and it is
pseudo-contractive, if we can choose $M = 1$ in
(\ref{semigroup-growth}).

Every $C_0$-semigroup is not far from being pseudo-contractive.
Indeed, for an arbitrary $s > 0$, we have, by
(\ref{semigroup-growth}), the estimate
\begin{align*}
\| S_t \| \leq e^{\omega(s) t}, \quad t \geq s
\end{align*}
where we have set $\omega(s) := \frac{\ln M}{s} + \omega$.
Nevertheless, there are $C_0$-semigroups, which are not
pseudo-contractive. For a counter example, we choose, following
\cite[Ex. I.5.7.iii]{Engel-Nagel}, the Hilbert space $H :=
L^2(\mathbb{R})$ and the shift semigroup $(S_t)_{t \geq 0}$ with
jump, defined as
\begin{align*}
S_t h(x) :=
\begin{cases}
2h(x+t), & x \in [-t,0]
\\ h(x+t), & \text{otherwise}
\end{cases}
\end{align*}
for $h \in H$. Then $(S_t)_{t \geq 0}$ is a $C_0$-semigroup on $H$
with $\| S_t \| = 2$ for all $t > 0$, because $\| S_t
\mathbbm{1}_{[0,t]} \| = 2 \| \mathbbm{1}_{[0,t]} \|$.

However, many semigroups of practical relevance are
pseudo-contractive, and then the following result shows that
Assumption \ref{ass-group} is satisfied.

\begin{proposition}\label{prop-group}
Assume the semigroup $(S_t)_{t \geq 0}$ is pseudo-contractive. Then
there exist another separable Hilbert space $\mathcal{H}$ and a $C_0$-group $(U_t)_{t \in \mathbb{R}}$ on $\mathcal{H}$ such
that (\ref{diagram-commutes}) is satisfied, where $\ell \in
L(H,\mathcal{H})$ is an isometric embedding and $\pi := \ell^* \in
L(\mathcal{H},H)$ is the orthogonal projection from $\mathcal{H}$
into $H$.
\end{proposition}

\begin{proof}
Since the semigroup $(S_t)_{t \geq 0}$ is pseudo-contractive, there
exists $\omega \geq 0$ such that (\ref{def-pseudo-contr}) is
satisfied. Hence, the $C_0$-semigroup $(T_t)_{t \geq 0}$ defined as
$T_t := e^{-\omega t} S_t$, $t \in \mathbb{R}_+$ is contractive. By
the Sz\H{o}kefalvi-Nagy theorem on unitary dilations (see e.g.
\cite[Thm. I.8.1]{Nagy}, or \cite[Sec. 7.2]{Davies}), there exist
another separable Hilbert space $\mathcal{H}$ and a unitary
$C_0$-group $(V_t)_{t \in \mathbb{R}}$ in $\mathcal{H}$ such that
\begin{align*}
\pi V_t \ell = T_t \quad \text{for all $t \in \mathbb{R}_+$,}
\end{align*}
where $\ell \in L(H,\mathcal{H})$ is an isometric
embedding and the adjoint operator $\pi := \ell^* \in
L(\mathcal{H},H)$ is the orthogonal projection from
$\mathcal{H}$ into $H$. Defining the $C_0$-group $(U_t)_{t
\in \mathbb{R}}$ as $U_t := e^{\omega t} V_t$, $t \in \mathbb{R}$
completes the proof.
\end{proof}

We suppose from now on Assumption \ref{ass-group}. There exist constants $M \geq 1$ and $\omega \in \mathbb{R}$ such
that
\begin{align}\label{group-growth}
\| U_t \| \leq M e^{\omega |t|}, \quad t \in \mathbb{R}
\end{align}
see \cite[p. 79]{Engel-Nagel}. Now let $\alpha : C_{\rm ad}(\mathbb{R}_+) \rightarrow H_{\mathcal{P}}$, $\sigma
: C_{\rm ad}(\mathbb{R}_+) \rightarrow (L_2^0)_{\mathcal{P}}$ and $\gamma :
C_{\rm ad}(\mathbb{R}_+) \rightarrow H_{\mathcal{P} \otimes \mathcal{E}}$ be given.
We suppose that Assumptions \ref{ass-0-strong}, \ref{ass-1-strong}, \ref{ass-2-strong} are satisfied.

In order to solve the stochastic
partial differential equation (\ref{general-SPDE}), we consider the
$\mathcal{H}$-valued stochastic differential equation
\begin{align}\label{equation-help}
\left\{
\begin{array}{rcl}
d R_t & = & \tilde{\alpha}(R)_t dt + \tilde{\sigma}(R)_t dW_t +
\int_E \tilde{\gamma}(R)(t,x) (\mu(dt,dx) - F(dx) dt)
\medskip
\\ R|_{[0,t_0]} & = & h,
\end{array}
\right.
\end{align}
where $t_0 \in \mathbb{R}_+$ and $h \in C_{\rm ad}([0,t_0];\mathcal{H})$, and where $\tilde{\alpha}
: C_{\rm ad}(\mathbb{R}_+;\mathcal{H}) \rightarrow \mathcal{H}_{\mathcal{P}}$,
$\tilde{\sigma} : C_{\rm ad}(\mathbb{R}_+;\mathcal{H}) \rightarrow
L_2(U_0,\mathcal{H})_{\mathcal{P}}$ and $\tilde{\gamma} : C_{\rm
ad}(\mathbb{R}_+;\mathcal{H}) \rightarrow \mathcal{H}_{\mathcal{P} \otimes
\mathcal{E}}$ are defined as
\begin{align}\label{def-alpha-new}
\tilde{\alpha}(R)_t &:= U_{-t}^{t_0} \ell \alpha(\pi U^{t_0} R)_t,
\\ \label{def-sigma-new} \tilde{\sigma}(R)_t &:= U_{-t}^{t_0} \ell \sigma(\pi U^{t_0} R)_t,
\\ \label{def-gamma-new} \tilde{\gamma}(R)(t,x) &:= U_{-t}^{t_0} \ell \gamma(\pi U^{t_0} R)(t,x).
\end{align}
In the above definitions, we have used the notation
\begin{align*}
U_t^{t_0} :=
\begin{cases}
U_{t-t_0}, & t \geq t_0
\\ {\rm Id}, & t \in (-t_0,t_0)
\\ U_{t_0+t}, & t \leq -t_0
\end{cases}
\end{align*}
and $\pi U^{t_0} R \in C_{\rm ad}(\mathbb{R}_+;H)$ denotes the process $(\pi U^{t_0} R)_t := \pi U_t^{t_0} R_t$, $t \geq 0$.
Note that $\tilde{\alpha}, \tilde{\sigma}, \tilde{\gamma}$ indeed
map into the respective spaces of predictable processes, because
$(t,h) \mapsto U_t h$ is continuous on $\mathbb{R} \times
\mathcal{H}$, see, e.g., \cite[Lemma VII.4.3]{Werner}. By (\ref{group-growth}), they also
fulfill Assumptions \ref{ass-0-strong}, \ref{ass-1-strong}, \ref{ass-2-strong}, where
the function $L$ is replaced by
\begin{align}\label{new-Lip-const}
L(t) \rightsquigarrow \| \ell \| \big( \mathbbm{1}_{[0,t_0)} + M^2 e^{2 \omega (t-t_0)} \mathbbm{1}_{[t_0,\infty)} \big) \| \pi \| L(t), \quad t \geq 0
\end{align}
According to Theorem \ref{thm-existence}, for each $h \in
\mathcal{C}_{\rm ad}([0,t_0];\mathcal{H})$ there exists a
unique solution $R \in \mathcal{C}_{\rm ad}(\mathbb{R}_+;\mathcal{H})$ for
(\ref{equation-help}) with c\`adl\`ag paths on $[t_0,\infty)$, and it satisfies
\begin{align}\label{solution-in-S2-larger}
\mathbb{E} \bigg[ \sup_{t \in [t_0,T]} \| R_t \|^2 \bigg] < \infty
\quad \text{for all $T \geq t_0$.}
\end{align}

\begin{theorem}\label{thm-mild-weak}
Suppose that Assumptions \ref{ass-0-strong}, \ref{ass-1-strong}, \ref{ass-2-strong} and \ref{ass-group}
are fulfilled. Then, for each $t_0 \in \mathbb{R}_+$ and $h \in
\mathcal{C}_{\rm ad}[0,t_0]$ there exists a unique
mild and weak solution $r \in \mathcal{C}_{\rm ad}(\mathbb{R}_+)$ for
(\ref{general-SPDE}) with c\`adl\`ag paths on $[t_0,\infty)$, and it satisfies
(\ref{solution-in-S2}). The solution is given by $r = \pi U^{t_0} R$, where $R
\in \mathcal{C}_{\rm ad}(\mathbb{R}_+;\mathcal{H})$ denotes the solution for
(\ref{equation-help}) with $R|_{[0,t_0]} = \ell h$.
\end{theorem}

\begin{proof}
Let $t_0 \in \mathbb{R}_+$ and $h \in
\mathcal{C}_{\rm ad}[0,t_0]$ be arbitrary.
The $H$-valued process $r := \pi U^{t_0} R$ belongs to $\mathcal{C}_{\rm ad}(\mathbb{R}_+)$, it
satisfies (\ref{solution-in-S2}) by virtue of
(\ref{solution-in-S2-larger}), and it has also c\`adl\`ag paths on $[t_0,\infty)$, because
$(t,h) \mapsto U_t h$ is continuous on $\mathbb{R}_+ \times
\mathcal{H}$, see, e.g., \cite[Lemma VII.4.3]{Werner}.
Using (\ref{diagram-commutes}) we obtain $r|_{[0,t_0]} = \pi \ell h = h$ and almost surely
\begin{align*}
r_t &= (\pi U^{t_0} R)_t = \pi U_{t-t_0} R_t 
\\ &= \pi U_{t-t_0} \bigg( \ell h_{t_0} + \int_{t_0}^t U_{t_0-s} \ell \alpha(\pi U^{t_0}
R)_s ds + \int_{t_0}^t U_{t_0-s} \ell \sigma(\pi U^{t_0} R)_s
dW_s
\\ &\quad + \int_{t_0}^t
\int_E U_{t_0-s} \ell \gamma(\pi U^{t_0} R)(s,x) (\mu(ds,dx) - F(dx) ds) \bigg)
\\ &= S_{t-t_0} h_{t_0} + \int_{t_0}^t S_{t-s} \alpha(r)_s ds + \int_{t_0}^t S_{t-s}
\sigma(r)_s dW_s
\\ &\quad + \int_{t_0}^t \int_E S_{t-s} \gamma(r)(s,x) (\mu(ds,dx) - F(dx)
ds), \quad t \geq t_0
\end{align*}
showing that $r$ is a mild solution for (\ref{general-SPDE}). By virtue of Lemma \ref{lemma-well-defined} we have
$\sigma(r) \in L^2(W;L_2^0)$ and $\gamma(r) \in L^2(\mu;H)$.
Applying Lemma \ref{lemma-mild-weak} proves that $r$ is also a weak
solution for (\ref{general-SPDE}).

For two mild solutions $r, \tilde{r} \in \mathcal{C}_{\rm ad}(\mathbb{R}_+)$
of (\ref{general-SPDE}), which are c\`adl\`ag on $[t_0,\infty)$, and an
arbitrary $T \geq t_0$, by using H\"older's inequality, the
It\^o-isometries (\ref{Ito-isometry-Wiener}),
(\ref{Ito-isometry-mu}) and the Lipschitz conditions
(\ref{Lip-alpha-SDE-strong}), (\ref{Lip-sigma-SDE-strong}),
(\ref{Lip-gamma-SDE-strong}), the inequality
\begin{align*}
\| r - \tilde{r} \|_{[t_0,t]}^2 &= \sup_{s \in [t_0,t]} \mathbb{E}[ \| r_s - \tilde{r}_s \|^2 ] 
\\ &\leq 3 M^2 e^{2 \omega (T-t_0)} (T-t_0
+ 2) \int_{t_0}^t L(v)^2 \| r - \tilde{r} \|_{[t_0,v]}^2 dv,
\quad t \in [t_0,T]
\end{align*}
is valid, where $M \geq 1$ and $\omega \in \mathbb{R}$ stem from
(\ref{group-growth}). Using the Gronwall Lemma and the hypothesis that
$r$ and $\tilde{r}$ are c\`adl\`ag on $[t_0,\infty)$, we conclude that $r$ and
$\tilde{r}$ are indistinguishable on $[t_0,\infty)$. Taking into account Lemma
\ref{lemma-weak-mild}, this proves the desired uniqueness of mild
and weak solutions for (\ref{general-SPDE}).
\end{proof}

\begin{remarks}\mbox{}
\begin{enumerate}
\item The idea to use the Sz\H{o}kefalvi-Nagy theorem on unitary
dilations in order to overcome the difficulties arising from
stochastic convolutions, is due to E.~Hausenblas and J.~Seidler, see
\cite{Seidler} and \cite{Seidler2}.

\item Imposing Assumptions \ref{ass-0-strong}, \ref{ass-1-Lp}, \ref{ass-2-Lp} and \ref{ass-group} we obtain the $L^p$-version of Theorem \ref{thm-mild-weak}.
\end{enumerate}
\end{remarks}

\begin{remark}
Another interpretation of Theorem \ref{thm-mild-weak} is the following: it is well known that generic mild (or weak) solutions of SPDEs \eqref{general-SPDE} are not Hilbert space valued semi-martingales due to lack of regularity in time of the finite variation part. However, our method shows that we can decompose every mild (or weak) solution as $ r_t = \pi U_{t-t_0} R_t $, where $ R $ is a semi-martingale, $ U $ a strongly continuous group and $ \pi $ the orthogonal projection due to Assumption \ref{ass-group}, and $ t \geq t_0 $.
\end{remark}

\section{Stability and regularity of stochastic partial differential
equations}\label{sec-stability-SPDE}

We shall now deal with stability and regularity of stochastic partial differential
equations of the kind (\ref{general-SPDE}). Stability and regularity
results for SPDEs can also be found in \cite{Ruediger-mild} and
\cite{Marinelli-Prevot-Roeckner}. Here, we can easily transfer the
results on stability from Section \ref{sec-SDE-stability} and on
regularity from Section \ref{sec-SDE-regularity} to SPDEs by the
method of the moving frame. For stability results, we provide the
details in this section.

As in Section \ref{sec-SPDE-ex}, we suppose Assumption \ref{ass-group} and that $\alpha : C_{\rm
ad}(\mathbb{R}_+) \rightarrow H_{\mathcal{P}}$, $\sigma : C_{\rm ad}(\mathbb{R}_+)
\rightarrow (L_2^0)_{\mathcal{P}}$ and $\gamma : C_{\rm ad}(\mathbb{R}_+)
\rightarrow H_{\mathcal{P} \otimes \mathcal{E}}$ fulfill Assumptions
\ref{ass-0-strong}, \ref{ass-1-strong}, \ref{ass-2-strong}. For each $n \in
\mathbb{N}$, let $\alpha_n : C_{\rm ad}(\mathbb{R}_+) \rightarrow
H_{\mathcal{P}}$, $\sigma_n : C_{\rm ad}(\mathbb{R}_+) \rightarrow
(L_2^0)_{\mathcal{P}}$ and $\gamma_n : C_{\rm ad}(\mathbb{R}_+) \rightarrow
H_{\mathcal{P} \otimes \mathcal{E}}$ be such that Assumptions
\ref{ass-3a-strong}, \ref{ass-3-strong}, \ref{ass-4-strong} are fulfilled.
Furthermore, let $t_0 \in \mathbb{R}_+$, $h \in \mathcal{C}_{\rm ad}[0,t_0]$
and for each $n \in \mathbb{N}$ let $h^n \in
\mathcal{C}_{\rm ad}[0,t_0]$ and $B_n \in \mathcal{E}$ be
given.

According to Theorem \ref{thm-mild-weak}, there exists a unique
solution $r \in \mathcal{C}_{\rm ad}(\mathbb{R}_+)$ for
(\ref{general-SPDE}) with $r|_{[0,t_0]} = h$ with c\`adl\`ag paths on $[t_0,\infty)$ satisfying
(\ref{solution-in-S2}), and for each $n \in \mathbb{N}$ there exists
a unique solution $r^n \in \mathcal{C}_{\rm ad}(\mathbb{R}_+)$ for
\begin{align*}
\left\{
\begin{array}{rcl}
dr_t^n & = & (Ar_t^n + \alpha(r^n)_t)dt + \sigma(r^n)_t dW_t
\\ && + \int_{B_n} \gamma(r^n)(t,x) (\mu(dt,dx) - F(dx) dt)
\medskip
\\ r^n|_{[0,t_0]} & = & h^n
\end{array}
\right.
\end{align*}
with c\`adl\`ag paths on $[t_0,\infty)$
satisfying $\mathbb{E} [ \sup_{t \in [t_0,T]} \| r_t^n \|^2 ] <
\infty$ for all $T \geq t_0$. We suppose that Assumption \ref{ass-convergence}, in which $r \in \mathcal{C}_{\rm ad}(\mathbb{R}_+)$ denotes the mild and weak solution for (\ref{general-SPDE}), holds true.

\begin{proposition}\label{prop-stability-mild}
Suppose that Assumptions \ref{ass-0-strong}, \ref{ass-1-strong}, \ref{ass-2-strong},
\ref{ass-3a-strong}, \ref{ass-3-strong}, \ref{ass-4-strong}, \ref{ass-convergence} and
\ref{ass-group} are fulfilled. Then, there exist maps $K_1,K_2 : \mathbb{R}_+ \rightarrow \mathbb{R}_+$, only depending on the Lipschitz function $L$, such that the following statements are valid:

\begin{enumerate}
\item If $h^n \rightarrow h$ in $C_{\rm ad}[0,t_0]$, then for each $T \geq t_0$ we have the estimate
\begin{align}\label{est-expectation-mild}
\sup_{t \in [0,T]} \mathbb{E}[ \| r_t - r_t^n \|^2 ] &\leq K_1 \big( \| h - h^n \|_{[0,t_0]}^2 + C_n^2 \big)
\rightarrow 0 \quad \text{for $n \rightarrow \infty$,}
\end{align}
where $K_1 = K_1(T)$ and $C_n = C_n(T,r)$ is defined in (\ref{def-Cn}).

\item If even $h^n \rightarrow h$ in $S^2[0,t_0]$, then for each $T \geq t_0$ we have the estimate
\begin{align}\label{est-sup-expectation-mild}
\mathbb{E}\bigg[ \sup_{t \in [0,T]} \| r_t - r_t^n \|^2 \bigg] &\leq
K_2 \big( \| h - h^n \|_{S^2[0,t_0]}^2 + C_n^2 \big)
\rightarrow 0 \quad \text{for $n \rightarrow \infty$,}
\end{align}
where $K_2 = K_2(T)$ and $C_n = C_n(T,r)$ is defined in (\ref{def-Cn}).
\end{enumerate}
\end{proposition}

\begin{proof}
We define
$\tilde{\alpha} : C_{\rm ad}(\mathbb{R}_+;\mathcal{H}) \rightarrow
\mathcal{H}_{\mathcal{P}}$, $\tilde{\sigma} : C_{\rm
ad}(\mathbb{R}_+;\mathcal{H}) \rightarrow L_2(U_0,\mathcal{H})_{\mathcal{P}}$ and
$\tilde{\gamma} : C_{\rm ad}(\mathbb{R}_+;\mathcal{H}) \rightarrow
\mathcal{H}_{\mathcal{P} \otimes \mathcal{E}}$ by
(\ref{def-alpha-new}), (\ref{def-sigma-new}), (\ref{def-gamma-new}).
Moreover, for each $n \in \mathbb{N}$, we define $\tilde{\alpha}_n :
C_{\rm ad}(\mathbb{R}_+;\mathcal{H}) \rightarrow \mathcal{H}_{\mathcal{P}}$,
$\tilde{\sigma}_n : C_{\rm ad}(\mathbb{R}_+;\mathcal{H}) \rightarrow
L_2(U_0,\mathcal{H})_{\mathcal{P}}$ and $\tilde{\gamma}_n : C_{\rm
ad}(\mathbb{R}_+;\mathcal{H}) \rightarrow \mathcal{H}_{\mathcal{P} \otimes
\mathcal{E}}$ as
\begin{align*}
\tilde{\alpha}_n(R)_t &:= U_{-t}^{t_0} \ell \alpha_n(\pi U^{t_0} R)_t,
\\ \tilde{\sigma}_n(R)_t &:= U_{-t}^{t_0} \ell \sigma_n(\pi U^{t_0} R)_t,
\\ \tilde{\gamma}_n(R)(t,x) &:= U_{-t}^{t_0} \ell \gamma_n(\pi U^{t_0} R)(t,x).
\end{align*}
According to Theorem \ref{thm-mild-weak}, we have $r = \pi U^{t_0} R$, where $R \in \mathcal{C}_{\rm ad}(\mathbb{R}_+;\mathcal{H})$ denotes the solution for
(\ref{equation-help}) with $R|_{[0,t_0]} = \ell h$, and for each $n \in
\mathbb{N}$ we have $r^n = \pi U^{t_0} R^n$, where $R^n \in \mathcal{C}_{\rm ad}(\mathbb{R}_+;\mathcal{H})$ denotes the solution for
\begin{align*}
\left\{
\begin{array}{rcl}
d R_t^n & = & \tilde{\alpha}(R^n)_t dt + \tilde{\sigma}(R^n)_t dW_t
+ \int_{B_n} \tilde{\gamma}(R^n)(t,x) (\mu(dt,dx) - F(dx) dt)
\medskip
\\ R^n|_{[0,t_0]} & = & \ell h^n.
\end{array}
\right.
\end{align*}
By (\ref{group-growth}), the coefficients
$\tilde{\alpha},\tilde{\sigma},\tilde{\gamma}$ and
$\tilde{\alpha}_n,\tilde{\sigma}_n,\tilde{\gamma}_n$, $n \in
\mathbb{N}$ fulfill Assumptions \ref{ass-0-strong}, \ref{ass-1-strong},
\ref{ass-2-strong}, \ref{ass-3a-strong}, \ref{ass-3-strong}, \ref{ass-4-strong}, where
the function $L$ is replaced by (\ref{new-Lip-const}). Moreover, by (\ref{group-growth}), for each $n \in \mathbb{N}$ we have
\begin{equation}\label{est-constants-C-n}
\begin{aligned}
\tilde{C}_n(T,R) &:= \bigg( \mathbb{E}\bigg[ \int_{t_0}^T \|
\tilde{\alpha}([R])_s - \tilde{\alpha}_n([R])_s \|^2 ds \bigg] 
\\ &\quad + \mathbb{E}\bigg[
\int_{t_0}^T \| \tilde{\sigma}([R])_s - \tilde{\sigma}_n([R])_s
\|_{L_2(U_0,\mathcal{H})}^2 ds \bigg]
\\ &\quad + \mathbb{E} \bigg[ \int_{t_0}^t \int_E \| \tilde{\gamma}([R])(s,x) -
\tilde{\gamma}_n([R])(s,x) \|^2 F(dx) ds \bigg] 
\\ &\quad +
\mathbb{E}\bigg[ \int_{t_0}^T \int_{E \setminus B_n} \|
\tilde{\gamma}([R])(s,x) \|^2 F(dx)ds \bigg] \bigg)^{\frac{1}{2}} \leq \| \ell \| M e^{\omega(T-t_0)} C_n(T,r)
\end{aligned}
\end{equation}
for all $T \geq t_0$.
In particular, Assumption \ref{ass-convergence} is fulfilled for $\tilde{\alpha}([R])$, $\tilde{\sigma}([R])$, $\tilde{\gamma}([R])$ and $\tilde{\alpha}_n([R])$, $\tilde{\sigma}_n([R])$, $\tilde{\gamma}_n([R])$, $n \in \mathbb{N}$.
If $h^n \rightarrow h$ in $C_{\rm ad}[0,t_0]$, then by applying Proposition \ref{prop-stability-strong} and noting (\ref{group-growth}) and (\ref{est-constants-C-n}), for each $T \geq t_0$ we obtain the estimate
\begin{align*}
&\sup_{t \in [0,T]} \mathbb{E}[ \| r_t - r_t^n \|^2 ] \leq \| \pi \|^2 M^2 e^{2 \omega(T-t_0)} \sup_{t \in [0,T]} \mathbb{E}[ \| R_t - R_t^n \|^2 ] 
\\ &\leq \| \pi \|^2 M^2 e^{2 \omega(T-t_0)} K_1(T) \big( \| h - h^n \|_{[0,t_0]}^2 + \tilde{C}_n(T,R)^2 \big) 
\\ &\leq \| \pi \|^2 M^2 e^{2 \omega(T-t_0)} K_1(T) \big( \| h - h^n \|_{[0,t_0]}^2 + \| \ell \|^2 M^2 e^{2 \omega(T-t_0)} C_n(T,r)^2 \big) \rightarrow 0
\end{align*}
for $n \rightarrow \infty$, where the map $K_1 : \mathbb{R}_+ \rightarrow \mathbb{R}_+$ stems from Proposition \ref{prop-stability-strong},
showing (\ref{est-expectation-mild}). Analogously, if $h^n \rightarrow h$ in $S^2[0,t_0]$, we get (\ref{est-sup-expectation-mild}).
\end{proof}

By Proposition \ref{prop-stability-mild}, the statement of Remark \ref{remark-solution-map} concerning the Lipschitz continuity of the solution map $h \mapsto r^h$ is also valid for SPDEs.

Analogously to stability results also the results on regularity can
be transferred to SPDEs by the method of the moving frame. The
arguments of Section \ref{sec-SDE-regularity} can be transferred
literally. The same arguments hold true for $ L^p $-estimates.

\section{Stochastic partial differential equations with state dependent coefficients}\label{sec-SPDE-markov}

In this section, we deal with stochastic partial differential
equations with state dependent coefficients, which may depend on the randomness $\omega$, the time $t$ and finitely many states of the path of the solution.
As we shall see, this is a special case
of the framework from Section \ref{sec-SPDE-ex}.

Let $K \in \mathbb{N}$ and $0 \leq \delta_1 < \ldots < \delta_K \leq 1$ be given. Moreover,
let $\alpha : \Omega \times \mathbb{R}_+ \times H^K \rightarrow H$, $\sigma :
\Omega \times \mathbb{R}_+ \times H^K \rightarrow L_2^0$ be $\mathcal{P} \otimes \mathcal{B}(H^K)$-measurable and $\gamma : \Omega \times \mathbb{R}_+
\times H^K \times E \rightarrow H$ be $\mathcal{P} \otimes \mathcal{B}(H^K) \otimes \mathcal{E}$-measurable.

\begin{assumption}\label{ass-1-several}
Denoting by $\mathbf{0} \in H^K$ the zero vector, we
assume that
\begin{align*}
t \mapsto \mathbb{E} [ \| \alpha(t,\mathbf{0}) \|^2] &\in \mathcal{L}_{\rm loc}^{1}(\mathbb{R}_+),
\\ t \mapsto \mathbb{E} [ \| \sigma(t,\mathbf{0}) \|_{L_2^0}^2] &\in \mathcal{L}_{\rm
loc}^{1}(\mathbb{R}_+),
\\ t \mapsto \mathbb{E} \bigg[ \int_E \| \gamma(t,\mathbf{0},x)
\|^2 F(dx) \bigg] &\in \mathcal{L}_{\rm loc}^{1}(\mathbb{R}_+).
\end{align*}
\end{assumption}

\begin{assumption}\label{ass-2-several}
We assume there is a function $L \in \mathcal{L}_{\rm loc}^{2}(\mathbb{R}_+)$
such that almost surely
\begin{align*}
\| \alpha(t,h_1) - \alpha(t,h_2) \| &\leq L(t) \sum_{i=1}^K \| h_1^i - h_2^i \|,
\\ \| \sigma(t,h_1) - \sigma(t,h_2) \|_{L_2^0} &\leq
L(t) \sum_{i=1}^K \| h_1^i - h_2^i \|,
\\ {\bigg (\int_E \| \gamma(t,h_1,x)
- \gamma(t,h_2,x) \|^2 F(dx) \bigg )}^{\frac{1}{2}} &\leq L(t) \sum_{i=1}^K \| h_1^i - h_2^i \|
\end{align*}
for all $t \in \mathbb{R}_+$ and all $h_1,h_2 \in H^K$.
\end{assumption}

\begin{corollary}\label{cor-several}
Suppose that Assumptions \ref{ass-group} and \ref{ass-1-several}, \ref{ass-2-several} are fulfilled. Then, for each $t_0 \in \mathbb{R}_+$ and $h \in
\mathcal{C}_{\rm ad}[0,t_0]$ there exists a unique mild and weak solution
$r \in \mathcal{C}_{\rm ad}(\mathbb{R}_+)$ for
\begin{align}\label{path-dep}
\left\{
\begin{array}{rcl}
dr_t & = & (Ar_t + \alpha(t,r_{\delta_1 t},\ldots,r_{\delta_K t}))dt + \sigma(t,r_{\delta_1 t},\ldots,r_{\delta_K t}) dW_t 
\\ && + \int_E
\gamma(t,r_{\delta_1 t-},\ldots,r_{\delta_K t-},x) (\mu(dt,dx) - F(dx) dt)
\medskip
\\ r|_{[0,t_0]} & = & h
\end{array}
\right.
\end{align}
with c\`{a}dl\`{a}g paths on $[t_0,\infty)$, and it satisfies (\ref{solution-in-S2}).
\end{corollary}

\begin{proof}
For every $r \in C_{\rm ad}(\mathbb{R}_+)$ let ${}^p r \in \mathcal{C}_{\rm ad}(\mathbb{R}_+)$ be a predictable representative of $r$, which, due to \cite[Prop. 3.6.ii]{Da_Prato}, always exists. Now, we define
the maps $\tilde{\alpha} : C_{\rm ad}(\mathbb{R}_+) \rightarrow
H_{\mathcal{P}}$, $\tilde{\sigma} : C_{\rm ad}(\mathbb{R}_+) \rightarrow
(L_2^0)_{\mathcal{P}}$ and $\tilde{\gamma} : C_{\rm ad}(\mathbb{R}_+)
\rightarrow H_{\mathcal{P} \otimes \mathcal{E}}$ by
\begin{align*}
\tilde{\alpha}(r)_t(\omega) &:= \alpha(\omega,t,{}^p r_{\delta_1 t}(\omega),\ldots,{}^p r_{\delta_K t}(\omega)), \quad (\omega,t) \in \Omega \times \mathbb{R}_+
\\ \tilde{\sigma}(r)_t(\omega) &:= \sigma(\omega,t,{}^p r_{\delta_1 t}(\omega),\ldots,{}^p r_{\delta_K t}(\omega)), \quad (\omega,t) \in \Omega \times \mathbb{R}_+
\\ \tilde{\gamma}(r)(t,x)(\omega) &:= \gamma(\omega,t,{}^p r_{\delta_1 t}(\omega),\ldots,{}^p r_{\delta_K t}(\omega),x), \quad (\omega,t,x) \in
\Omega \times \mathbb{R}_+ \times E.
\end{align*}
Note that for a predictable process $(r_t)_{t \geq 0}$ and an arbitrary $0 \leq \delta \leq 1$ the process $(r_{\delta t})_{t \geq 0}$ is predictable, too. Hence,
$\tilde{\alpha}$, $\tilde{\sigma}$, $\tilde{\gamma}$ indeed map into the respective spaces of predictable processes, because $\alpha$, $\sigma$ are $\mathcal{P} \otimes \mathcal{B}(H^K)$-measurable and $\gamma$ is $\mathcal{P} \otimes \mathcal{B}(H^K) \otimes \mathcal{E}$-measurable. Assumption \ref{ass-0-strong} holds true by the definition of $\tilde{\alpha}$, $\tilde{\sigma}$, $\tilde{\gamma}$ and Assumption \ref{ass-1-strong} is satisfied by Assumption \ref{ass-1-several}. Using Assumption \ref{ass-2-several}, for all $t \in \mathbb{R}_+$ and all $r^1,r^2 \in C_{\rm ad}(\mathbb{R}_+)$ we obtain
\begin{align*}
&\mathbb{E}[\| \alpha(r^1)_t - \alpha(r^2)_t \|^2] = \mathbb{E} [\| \alpha(t,r_{\delta_1 t}^1,\ldots,r_{\delta_K t}^1) - \alpha(t,r_{\delta_1 t}^2,\ldots,r_{\delta_K t}^2) \|^2]
\\ &\leq L(t)^2 \mathbb{E} \Bigg[ \bigg( \sum_{i=1}^K \| r_{\delta_i t}^1 - r_{\delta_i t}^2 \| \bigg)^2 \Bigg] \leq K L(t)^2 \sum_{i=1}^K \mathbb{E}[ \| r_{\delta_i t}^1 - r_{\delta_i t}^2 \|^2 ]
\\ &\leq K^2 L(t)^2 \max_{i=1,\ldots,K} \mathbb{E}[\| r_{\delta_i t}^1 - r_{\delta_i t}^2 \|^2] \leq K^2 L(t)^2 \| r^1 - r^2 \|_{[0,t]}^2.
\end{align*}
An analogous argumentation for $\sigma$ and $\gamma$ proves that Assumption \ref{ass-2-strong} is fulfilled. Applying Theorem \ref{thm-mild-weak}, there exists a unique mild and weak solution
$r \in \mathcal{C}_{\rm ad}(\mathbb{R}_+)$ for (\ref{general-SPDE}) with c\`{a}dl\`{a}g paths on $[t_0,\infty)$ satisfying (\ref{solution-in-S2}). For every $T \geq t_0$ we have, by using Assumption \ref{ass-2-several}, and since each path of $r$ has only countably many jumps on the interval $[t_0,T]$,
\begin{align*}
&\mathbb{E} \bigg[ \int_{t_0}^T \int_E \| \gamma(t,{}^p r_{\delta_1 t},\ldots,{}^p r_{\delta_K t},x) - \gamma(t,r_{\delta_1 t-},\ldots,r_{\delta_K t-},x) \|^2 F(dx) dt \bigg]
\\ &\leq \int_{t_0}^T \mathbb{E} \Bigg[ L(t)^2 \bigg( \sum_{i=1}^K \| {}^p r_{\delta_i t} - r_{\delta_i t-} \| \bigg)^2 \Bigg] dt \leq K \sum_{i=1}^K \int_{t_0}^T \mathbb{E} \big[ L(t)^2 \| r_{\delta_i t} - r_{\delta_i t-} \|^2 \big] dt
\\ &= K \sum_{i=1}^K \mathbb{E} \bigg[ \int_{t_0}^T L(t)^2 \| \Delta r_{\delta_i t} \|^2 dt \bigg] = 0.
\end{align*}
Therefore, $\gamma(t,{}^p r_{\delta_1 t},\ldots,{}^p r_{\delta_K t},x) \mathbbm{1}_{[t_0,\infty)}$ and $\gamma(t,r_{\delta_1 t-},\ldots,r_{\delta_K t-},x) \mathbbm{1}_{[t_0,\infty)}$ coincide in the space $L^2(\mu;H)$. Consequently, the process $r$ is also the unique mild and weak solution for (\ref{path-dep}).
\end{proof}

As a particular case, we now turn to the Markovian framework. Let $\alpha : \mathbb{R}_+ \times H \rightarrow H$, $\sigma :
\mathbb{R}_+ \times H \rightarrow L_2^0$ and $\gamma : \mathbb{R}_+
\times H \times E \rightarrow H$ be measurable.

\begin{assumption}\label{ass-1-mild}
We assume that
\begin{align*}
t \mapsto \| \alpha(t,0) \| &\in \mathcal{L}_{\rm loc}^{2}(\mathbb{R}_+),
\\ t \mapsto \| \sigma(t,0) \|_{L_2^0} &\in \mathcal{L}_{\rm
loc}^{2}(\mathbb{R}_+),
\\ t \mapsto \int_E \| \gamma(t,0,x)
\|^2 F(dx) &\in \mathcal{L}_{\rm loc}^{1}(\mathbb{R}_+).
\end{align*}
\end{assumption}

\begin{assumption}\label{ass-2-mild}
We assume there is a function $L \in \mathcal{L}_{\rm loc}^{2}(\mathbb{R}_+)$
such that
\begin{align*}
\| \alpha(t,h_1) - \alpha(t,h_2) \| &\leq L(t) \| h_1 - h_2 \|,
\\ \| \sigma(t,h_1) - \sigma(t,h_2) \|_{L_2^0} &\leq
L(t) \| h_1 - h_2 \|,
\\ {\bigg (\int_E \| \gamma(t,h_1,x)
- \gamma(t,h_2,x) \|^2 F(dx) \bigg )}^{\frac{1}{2}} &\leq L(t) \| h_1 - h_2 \|
\end{align*}
for all $t \in \mathbb{R}_+$ and all $h_1, h_2 \in H$.
\end{assumption}

\begin{corollary}\label{cor-SPDE}
Suppose that Assumptions \ref{ass-group} and \ref{ass-1-mild}, \ref{ass-2-mild} are fulfilled. Then, for each $h_0 \in
\mathcal{L}^2(\Omega,\mathcal{F}_0,\mathbb{P};H)$ there exists a unique
c\`adl\`ag, adapted, mean-square continuous mild and weak solution
$(r_t)_{t \geq 0}$ for
\begin{align*}
\left\{
\begin{array}{rcl}
dr_t & = & (Ar_t + \alpha(t,r_t))dt + \sigma(t,r_t) dW_t + \int_E
\gamma(t,r_{t-},x) (\mu(dt,dx) - F(dx) dt)
\medskip
\\ r_0 & = & h_0,
\end{array}
\right.
\end{align*}
and it satisfies (\ref{solution-in-S2}).
\end{corollary}

\begin{proof}
The assertion follows from Corollary \ref{cor-several} with $K = 1$, $\delta_1 = 1$ and $t_0 = 0$.
\end{proof}

We close this section with the time-homogeneous case. Let $\alpha :
H \rightarrow H$, $\sigma : H \rightarrow L_2^0$ and $\gamma : H
\times E \rightarrow H$ be measurable.

\begin{assumption}\label{ass-1-mild-hom}
We assume $\int_E \| \gamma(0,x) \|^2 F(dx) < \infty$.
\end{assumption}

\begin{assumption}\label{ass-2-mild-hom}
We assume that there is a constant $L \geq 0$ such that
\begin{align*}
\| \alpha(h_1) - \alpha(h_2) \| &\leq L \| h_1 - h_2 \|,
\\ \| \sigma(h_1) - \sigma(h_2) \|_{L_2^0} &\leq
L \| h_1 - h_2 \|,
\\ { \bigg ( \int_E \| \gamma(h_1,x)
- \gamma(h_2,x) \|^2 F(dx) \bigg )}^{\frac{1}{2}} &\leq L \| h_1 - h_2 \|
\end{align*}
for all $h_1, h_2 \in H$.
\end{assumption}

\begin{corollary}\label{cor-SPDE-hom}
Suppose that Assumptions \ref{ass-group} and \ref{ass-1-mild-hom}, \ref{ass-2-mild-hom}
are fulfilled. Then, for each $h_0 \in
\mathcal{L}^2(\Omega,\mathcal{F}_0,\mathbb{P};H)$ there exists a unique
c\`adl\`ag, adapted, mean-square continuous mild and weak solution
$(r_t)_{t \geq 0}$ for
\begin{align*}
\left\{
\begin{array}{rcl}
dr_t & = & (Ar_t + \alpha(r_t))dt + \sigma(r_t) dW_t + \int_E
\gamma(r_{t-},x) (\mu(dt,dx) - F(dx) dt)
\medskip
\\ r_0 & = & h_0,
\end{array}
\right.
\end{align*}
and it satisfies (\ref{solution-in-S2}).
\end{corollary}

\begin{proof}
The statement is an immediate consequence of Corollary \ref{cor-SPDE}.
\end{proof}

\begin{remark}\label{remark-Lp-state}
An analogous reasoning also provides the corresponding $L^p$-versions of
Corollaries \ref{cor-several}, \ref{cor-SPDE} and \ref{cor-SPDE-hom}.
\end{remark}

\begin{remark}
The time-inhomogeneous case can be considered by an extension of the
state space from $H$ to $ \mathbb{R} \times H $. However, one has to
pay attention at the boundary points of the interval $ [0,T] $,
where the vector fields have to be extended to the whole real line.
Nevertheless we shall consider in the setting of our numerical
applications the time-homogeneous case as the most characteristic one
for all further applications.
\end{remark}

\section{High-order (explicit-implicit) numerical schemes for stochastic partial differential equations with weak convergence order}\label{sec-SPDE-numerics}

We sketch in this last section high-order explicit-implicit numerical schemes for
stochastic partial differential equations with state-dependent coefficients as introduced in Section \ref{sec-SPDE-markov}. In this section (and only here) we will actually use that the Wiener process and the Poisson random measure are independent, see Section \ref{sec-independence}. By the stability results from Section \ref{sec-stability-SPDE} we can reduce the
problem to simpler driving signals, namely a finitely active Poisson random measure and to a finite number of
driving Wiener processes. We apply the results of \cite{BayTei:08} for the
time-dependent SDE, which -- due to the ``method of the moving
frame'' -- can be transferred to the general SPDE case. Our main
focus here is to work out so-called cubature schemes, extended by finite activity jump parts,
for time-dependent SDEs and therefore -- by the method of the moving frame -- for SPDEs. This 
also allows for high order numerical approximation schemes. Notice that cubature schemes are very adapted to SPDEs, since 
every local step can -- in contrast to Taylor schemes -- preserve the regularity of the states. Additionally a simple complexity analysis in the case of SPDEs yields that having a small amount $ p $ of local high order steps is cheaper than having a large amount $ p $ of local low order steps. The reason is that every local step means practically to solve a PDE numerically.

For this purpose we apply the respective notions from Section \ref{sec-SPDE-def}, \ref{sec-SPDE-ex} and \ref{sec-SPDE-markov} in order to formulate our conditions on the vector fields. Having the stability results of Section \ref{sec-stability-SPDE} in mind we do assume finite activity of the Poisson random measure, i.e.~$ F(E) < \infty $ and a finite dimensional Wiener process. Notice that this also allows for a statement on the rate of convergence to the original equations with possibly infinitely active jumps and infinitely many Brownian motions.

We consider here SPDEs of the type
\begin{align}\label{spde-numerics}
\left\{
\begin{array}{rcl}
dr_t & = & (Ar_t + \alpha(r_t))dt + \sigma(r_t) dW_t + \int_E
\gamma(r_{t-},x) (\mu(dt,dx) - F(dx) dt)
\medskip
\\ r_0 & \in & H
\end{array}
\right.
\end{align}

Let $ T > 0 $ denote a time-horizon. As in Section \ref{sec-SPDE-markov} we introduce measurable maps $\alpha :  H \rightarrow H$, $\sigma :  H \rightarrow L_2^0$ and $\gamma :  H \times E \rightarrow H$ and define maps $\tilde{\alpha} : [0,T] \times \mathcal{H} \rightarrow \mathcal{H}$, $\tilde{\sigma} : [0,T] \times \mathcal{H} \rightarrow \mathcal{H} $ and $\tilde{\gamma} : [0,T] \times \mathcal{H} \times E \rightarrow \mathcal{H} $ defined as
\begin{align*}
\tilde{\alpha}(t,R) &:= U_{-t} \ell \alpha(\pi U_t R), \quad t \in [0,T]
\\ \tilde{\sigma}(t,R) &:= U_{-t} \ell \sigma(\pi U_t R), \quad t \in [0,T]
\\ \tilde{\gamma}(t,R,x) &:= U_{-t} \ell \gamma(\pi U_t R,x), \quad (t,x) \in
[0,T] \times E.
\end{align*}

\begin{assumption}\label{ass-euler}
Fix $ m \geq 2 $ (a degree of accuracy for the high order scheme) and $ T > 0 $. We assume that for the vector fields $ \tilde{\alpha} $, $ \tilde{\sigma} $, $ \tilde{\gamma} $ there is a constant $ M > 0 $ such that for every radius $ C_1 > 0 $ we have
\begin{align*}
 \sup_{t \in [0,T], \, h \in \mathcal{H}, \| h \| \leq C_1} \| \partial^{k_1}_t \partial^{k_2}_h \tilde{\alpha}(t,  h) \| & \leq M C_1^{k_1},
\\ \sup_{t \in [0,T], \, h \in \mathcal{H}, \|h \| \leq C_1} \| \partial^{k_1}_t \partial^{k_2}_h \tilde{\sigma}(t, h) \|_{L_2^0} & \leq M C_1^{k_1},
\\ \sup_{t \in [0,T], \, h \in \mathcal{H}, \| h \| \leq C_1} \int_E \| \partial^{k_1}_t \partial^{k_2}_h \tilde{\gamma}(t,  h,x) \|^2 F(dx) &   \leq M C_1^{k_1},
\end{align*}
holds true for all $ k_1 + k_2 \leq m+1 $. In words, the growth of derivatives of the time-dependent vector fields up to order $  m+ 1 $ is polynomial in the radius of order $ k_1$, when $ k_1 $ denotes the order of the time derivative.

\end{assumption}

\begin{example}
Typical examples for vector fields satisfying the Assumptions \ref{ass-euler} are those of functional form (as applied in \cite{filipoteich}), i.e., choose a smooth map with all derivatives bounded $ \phi: \mathbb{R}^n \to D(\mathcal{A}^{\infty}) $, where $ \mathcal{A} $ denotes the infinitesimal generator of $ U $, and choose $ \xi_1,\cdots,\xi_n \in D((\mathcal{A^{*}})^{\infty}) $, then
\begin{equation}
\tilde{\sigma}(t,h)=  U_{-t} \ell \phi(\langle \xi_1 , \pi U_t h \rangle, \cdots,\langle \xi_n , \pi U_t h \rangle ) \label{functional-form}
\end{equation}
for $ h \in \mathcal{H} $ satisfies Assumptions \ref{ass-euler} for every $ m \geq 2 $.
\end{example}

\begin{remark}
As outlined in \cite{Tei:08} the previous assumptions imply $ \operatorname{Lip}(m+1) $-conditions on the ball with radius $ C_1 $ for the corresponding time-dependent vector fields $ \tilde{\alpha}, \tilde{\sigma}, \tilde{\gamma} $. This has an important meaning for the extension of our theory towards rough paths, see \cite{Tei:08}.
\end{remark}

\begin{remark}
The Assumptions \ref{ass-euler} lead to global existence and uniqueness of the corresponding time-dependent stochastic differential equations and the corresponding SPDEs \eqref{spde-numerics}. Due to the finite activity of the jump process the conditions also lead to existence of moments of any order of the solution process.
\end{remark}

We do first assume $ \gamma=0$, such that we find ourselves in a pure diffusion case. Furthermore we have assumed that the driving Wiener noise is finite dimensional, in other words we can write the stochastic partial differential equation in the moving frame and on the extended phase space $ \mathbb{R} \times \mathcal{H} $:
\begin{align}\label{time-dependent-system}
dR_t = \tilde{\alpha}(s,R_t) dt + \sum_{i=1}^d \tilde{\sigma}_i (s,R_t) dW^i_t, \, ds = dt,\\
R_0 = r_0, \, s_0=0.
\end{align}

Let us fix $ m \geq 2 $. Notice that the Assumption \ref{ass-euler} implies Assumptions \ref{ass-1-strong} and \ref{ass-2-strong}, in particular the vector fields are $(m+1)$-times differentiable in all variables. This allows us to state the standard result on short-time asymptotic for the stochastic differential equation \eqref{time-dependent-system}.

We now apply the notations for iterated stochastic integrals, i.e., the abbreviation
$$
W^{(i_1,\ldots,i_k)}_t = \int_{0 \leq t_1 \leq \dots \leq t_k \leq t} \circ d W^{i_1}_{t_1} \dots \circ d W^{i_k}_{t_k}
$$
is short for iterated Stratonovich integrals, where we apply $ \circ dW^0_t = dt $. Notice also the degree mapping
$$
\deg(i_1,\ldots,i_k) = k + \operatorname{card} \{ j \, | \; i_j =0  \},
$$
which counts any appearance of $ 0 $ in the multiindex $ (i_1,\ldots,i_k) $ twice, since $ dt = \circ dW^0_t = dt $ comes with twice the order of short time asymptotics than a Brownian motion. This is also the reason for the particular structure of the Assumptions \ref{ass-euler}. Recall also that any vector field $\sigma$ can be interpreted as a first order differential operator on test functions $f$ by
\begin{equation*}
  (\sigma f)(x) = Df(x)\cdot\sigma(x),\quad x\in H,
\end{equation*}
which will be applied extensively in the sequel.
\begin{theorem}\label{stoch_taylor_expansion}
Let $ g :\mathbb{R} \times \mathcal{H} \to \mathbb{R} $ be a smooth function with all derivatives bounded. Then we have the following asymptotic formula,
\begin{equation}
    g(t,R_t) = \sum_{\deg(i_1,\ldots,i_k) \leq m} (\beta_{i_1}\cdots \beta_{i_k}g)(0,R_0)
    W^{(i_1,\ldots,i_k)}_t + R_m(t,g,R_0), \quad R_0 \in \mathcal{H},
  \end{equation}
  with
  \begin{align}
    \sqrt{E(R_m(t,g,R_0)^2)} &\leq C t^{\frac{m+1}{2}}
    \max_{m<\deg(i_1,\ldots,i_k)\leq m+2} \; \sup_{0 \leq s \leq t}
    \sqrt{E(\abs{\beta_{i_1}\cdots\beta_{i_k}g(R_s,s)}^2)} \leq \\
& \leq \tilde{M} \sup_{2k_1 \leq m+1} \sup_{ 0 \leq s \leq t} E( {|| R_s ||}^{k_1}) < \infty,
  \end{align}
where $ \tilde{M} $ is a constant derived from Assumptions \eqref{ass-euler}.
\end{theorem}

\begin{proof}
The proof is a direct consequence of the results of \cite[Prop. 3.1]{BayTei:08}, where one additionally observes the necessary degrees of differentiability which are needed for the result. Notice in particular that the remainder term stays bounded due to the conditions of Assumption \ref{ass-euler}, in particular due to polynomial growth of the derivatives of at most order $ m + 1 $ and the existence of moments up to order $ m + 1 $ of the solution process.
\end{proof}

\begin{example} We formulate the short-time asymptotic formula in the case $m=2$ by taking the definitions of the vector fields $ \beta_0,\ldots,\beta_d $ and a smooth test function $ g: \mathcal{H} \rightarrow \mathbb{R} $, which does not depend on the additional (time-)state $ s $, then
\begin{align*}
g(R_t) = &  g(R_0) + \beta_0 g(R_0) t + \sum_{i=1}^d \beta_i g(R_0) W^i_t + \\
& \sum_{i,j=1}^d \beta_i \beta_j g(R_0) W^{(i,j)}_t  + \mathcal{O}(t^{\frac{3}{2}})\\
 = & g( R_0) + Dg(R_0)\bullet \alpha(R_0) t + \sum_{i=1}^d Dg(R_0) \bullet \sigma_i(R_0) W^i_t + \\
& + \sum_{i=1}^d Dg(R_0) \bullet (D \sigma_i(R_0) \bullet \sigma_i(R_0))\frac{{(W^i_t)}^2 - t}{2} + \\
& + \sum_{i \neq j=1}^d Dg(R_0) \bullet ( D \sigma_i(R_0) \bullet \sigma_j(R_0)) W^{(i,j)}_t + \mathcal{O}(t^{\frac{3}{2}})
\end{align*}
for $ t \geq 0 $, $ R_0 \in \mathcal{H} $. Heading for a strong Euler-Maruyama-scheme the previous formula yields -- formally evaluated for $ g = id $ -- the first iteration step from $ 0 \to t $
$$
R_0 \mapsto R_0 + \alpha(R_0) t + \sum_{i=1}^d \sigma_i(R_0) W^i_t.
$$ 
For the next step in the iteration we need the asymptotic expansion at time $ t $ and therefore also the vector fields $ \tilde{\alpha}, \tilde{\sigma} $ appear at time $ t $, namely
$$
R_t \mapsto R_t + \tilde{\alpha}(t,R_t) t + \sum_{i=1}^d \tilde{\sigma}_i(t,R_t) W^i_t.
$$
However, when one transfer the iteration of these two steps via $ \pi U_{2t} $ to $ H $ the described cancellation happens and one obtains the two-fold iteration of the time-homogeneous scheme
$$
r_0 \mapsto S_t r_0 + S_t \alpha(r_0) t + \sum_{i=1}^d S_t \sigma_i(r_0) W^i_t.
$$
Notice that this scheme is implicit in the linear PDE-part and
explicit in the stochastic components and the non-linear drift
component. Notice also that the weak convergence order $ 1 $ is
obtained if the Assumptions \ref{ass-euler} for $ m = 2 $ are
satisfied for smooth test functions with all derivatives bounded.
\end{example}

As explained in the literature, for instance in \cite{BayTei:08} or \cite{kloe/pla92}, we can derive high-order schemes (strong or weak) from the given short time-asymptotic expansion. The weak order of convergence -- given a short-time asymptotics of order $ t^{\frac{m+1}{2}} $ -- is then $ \frac{m-1}{2} $. Therefore we obtain high-order Taylor schemes for the time-dependent system \eqref{time-dependent-system}. However, even though possible, those Taylor schemes are usually not interesting -- except for the case $ m=2 $ -- since one has to work at each step with time derivatives of the vector fields, which corresponds to working with the infinitesimal generator of the semigroup.
\begin{example}
Consider a vector field $ \tilde{\sigma} $ of functional form \eqref{functional-form}, then apparently the time-derivative of the vector field, which appears in the stochastic Taylor expansion for $ m \geq 3 $, has the formula
\[
\frac{\partial}{\partial t} \tilde{\sigma}(t,R) = D \sigma(\langle \xi_1 , \pi U_t \rangle, \cdots,\langle \xi_n , \pi U_t \rangle) \bullet
(\langle \pi \mathcal{A}^{*} \xi_1 , U_t  h\rangle, \cdots,\langle \pi \mathcal{A}^{*} \xi_n ,  U_t h\rangle  ),
\]
which contains the infinitesimal generator and which is linearly growing in $ h $.
\end{example}

We present here a method to circumvent the problem that in each local step the infinitesimal generator appears, namely the cubature method: its implementation and structure work in the case of Hilbert space valued SDEs of type \eqref{time-dependent-system} in precisely the same way as in the finite dimensional case (see for instance \cite{BayTei:08} for details), since we do not have to deal with the unbounded infinitesimal generator. Convergence of global order $ \frac{m-1}{2} $ follows from Assumptions \ref{ass-euler} on any bounded set. On the other hand, each local time step does not contain derivatives of the vector fields in question, and preserves therefore the regularity of the state vector. We need one analytical preparation for this, namely the following lemma which tells that -- under Assumption \ref{ass-euler} -- we can suppose that on each bounded set there is a $ \Lip(m+1)$-extension of the vector fields on the whole extended phase space.

\begin{lemma}\label{lip-extension}
Define vector fields $ \beta_i $ on the extended phase space $ [0,T] \times \mathcal{H} $ by the following formulas:
\begin{align}
\beta_0(s,R) = \bigl(1, \tilde{\alpha}(s,R) - \frac{1}{2} \sum_{i=1}^d D \tilde{\sigma}_i(s,R) \bullet \tilde{\sigma}_i(s,R) \bigr), \\
\beta_i (s,R) = \bigl(0,  \tilde{\sigma}_i (s,R) \bigr),
\end{align}
for $ i =0,\ldots,d $. Then for each $ C_1 $ we find vector fields $ \beta_0^{C_1},\ldots,\beta_d^{C_1} $ which conincide with the previous vector fields on the ball with radius $ C_1 $ but are $ \Lip(m+1) $ on the whole extended phase space $ \mathbb{R} \times \mathcal{H} $.
\end{lemma}
\begin{proof}
This is not a consequence of the hard Whitney extension theorem but simply due to the fact that on $ C_1 + 1 $ we the vector fields $ \beta_i $ also satisfy a $ \Lip(m+1) $ condition. Multiplying with a bump function being equal to one on the ball of radius $ C_1 $ and vanishing outside radius $ C_1 + 1 $ yields the result.
\end{proof}

In the rest of the section we develop the necessary terminology for cubature methods: Theorem~\ref{stoch_taylor_expansion} shows that iterated Stratonovich integrals play the same r\^{o}le as polynomials play in deterministic Taylor expansion. Consequently, it is natural to use them in order to define cubature formulas. Let $C_{bv}([0,t];\R^d)$ denote the space of continuous paths of bounded
variation taking values in $\R^d$. As for the Brownian motion, we append a
component $\omega^0(t)=t$ for any $\omega\in C_{bv}([0,t];\R^d)$. Furthermore,
we establish the following convention: whenever $r_t$ is the solution to
some stochastic differential equation driven by Brownian motions $W$,
whether on a finite or infinite dimensional space, and $\omega\in
C_{bv}([0,t];\R^d)$, we denote by $r_t(\omega)$ the solution of the
deterministic differential equation given by formally replacing all
occurrences of ``$\circ dW^i_s$'' with ``$d\omega^i(s)$'' (with the same
initial values). Note that it is necessary that the SDE for $r$ is formulated in the Stratonovich sense (recall that the Stratonovich formulation does not necessarily make sense). With the following simple lemma we see
that time dependent coordinate transforms commute with the procedure of
replacing Brownian motions by deterministic trajectories.

\begin{lemma}\label{transfer-sde-spde}
Let $ \omega: [0,T] \to \mathbb{R}^d $ be a continuous curve with finite total variation. Then the time-dependent ordinary differential equation
\begin{align}
dR_t(\omega) = (\tilde{\alpha}(s,R_t(\omega)) - \frac{1}{2} \sum_{i=1}^d D \tilde{\sigma}_i(R_t(\omega),t) \bullet \tilde{\sigma}_i(t,R_t(\omega)))dt + \\
+ \sum_{i=1}^d \tilde{\sigma}_i (t,R_t(\omega)) d \omega^i(t),\, R_0 = r_0,
\end{align}
driven by $ \omega $ instead of the finite dimensional Wiener process $ W $, has a strong solution, which transfers via $ r_t(\omega) = \pi U_t R_t(\omega) $ to a mild solution of
\begin{equation} 
dr_t(\omega) = (A r_t(\omega) + \alpha(r_t(\omega)) - \frac{1}{2} \sum_{i=1}^d D \sigma_i(r_t(\omega)) \bullet \sigma_i(r_t(\omega)))dt + \sum_{i=1}^d \sigma_i (r_t(\omega)) d \omega^i(t).
\end{equation} 
\end{lemma}

Having in mind that one replaces Brownian motion $ W $ by a finite set of deterministic curves appearing with certain probabilities, we have to keep track of necessary moment conditions for (high-order) weak convergence, which is done in the following definition:

\begin{definition}
  \label{def:cubature_wiener_space}
  Fix $t>0$ and $m\geq1$. Positive weights $\lambda_1,\ldots,\lambda_N$
  summing up to $1$ and paths $\omega_1,\ldots,\omega_N\in C_{bv}([0,t];\R^d)$
  form a \emph{cubature formula on Wiener space} of degree $m$ if for all
  multi-indices $(i_1,\ldots,i_k)\in \mathcal{A}$ with
  $\deg(i_1,\ldots,i_k)\leq m$, $k\in\N$, we have that
  \begin{equation*}
    E(W^{(i_1,\ldots,i_k)}_t) = \sum_{l=1}^N \lambda_l
    W^{(i_1,\ldots,i_k)}_t(\omega_l), 
  \end{equation*}
  where we used the convention in line with the previous one, namely
$$
W^{(i_1,\ldots,i_k)}_t(\omega)= \int_{0 \leq t_1 \leq \cdots \leq t_k \leq t}
d \omega^{i_1}(t_1) \cdots d \omega^{i_k}(t_k).
$$ 
\end{definition}

Lyons and Victoir \cite{lyo/vic04} show the existence of cubature formulas on Wiener space for any $d$ and size $N \leq \#\{I\in\mathcal{A} | \deg(I)\leq m\}$ by applying Chakalov's theorem on cubature formulas and Chow's theorem
for nilpotent Lie groups. Moreover, due to the scaling properties of Brownian
motion (and its iterated Stratonovich integrals), i.e.
\begin{equation*}
  W^{(i_1,\ldots,i_k)}_t =^{\text{law}} \sqrt{t}^{\deg(i_1,\ldots,i_k)}
  W^{(i_1,\ldots,i_k)}_1, 
\end{equation*}
it is sufficient to construct cubature paths for $t=1$.

\begin{assumption}
  \label{ass:one-cubature-formula}
  Once and for all, we fix one cubature formula $\widetilde{\omega}_1, \ldots,
  \widetilde{\omega}_N$  with weights  $\lambda_1, \ldots, \lambda_N$ of degree
  $m \geq 2$ on the interval $[0,1]$. Without loss of generality we assume that $ \widetilde{\omega_i}(0)=0 $. 
  By abuse of notation, for any $t>0$, we will denote $\omega_l(s) = \sqrt{t} \widetilde{\omega}_l(s/t)$, $s \in [0,t]$,
  $l = 1,\ldots, N$, which yields a cubature formula for $[0,t]$ with the same weights $ \lambda_1, \ldots, \lambda_N $.
\end{assumption}

\begin{example}\label{ex:cubature_paths}
  For $d=1$ Brownian motions, a cubature formula on Wiener space of degree
  $m=3$ is given by $N=2$ paths
  \begin{equation*}
    \omega_1(s) = -\f{s}{\sqrt{t}},\ \omega_2(s) = \f{s}{\sqrt{t}}
  \end{equation*}
  for fixed time horizon $t$. The corresponding weights are given by
  $\lambda_1 = \lambda_2 = \f{1}{2}$.
\end{example}

When we deal with $ \Lip(m+1) $ vector fields on extended phase space we can write down -- by means of the finitely many cubature trajectories -- a local scheme. Notice that we have to replace the original vector fields $ \beta_0,\ldots,\beta_d $ by globally $ \Lip(m+1) $ vector fields $ \beta_0^{C_1},  \ldots, \beta_d^{C_1} $ on some large ball of radius $ C_1 $, see Lemma \ref{lip-extension}. The respective solutions of the SDEs are denoted by $ R^{C_1} $. Combining then the stochastic Taylor expansion, the deterministic Taylor expansion for solutions of ODEs driven by $ \omega_i $ for a cubature formula on Wiener space one obtains a one-step scheme for weak approximation of equations of type \eqref{time-dependent-system} precisely the same way as in \cite{lyo/vic04}. Indeed, we get
\begin{multline}\label{eq:bounded_1step}
\sup_{r_0 \in \mathcal{H}}\abs{E(g(t,R^{C_1}_t)) - \sum\nolimits_{l=1}^N \lambda_l g(t,R^{C_1}_t(\omega_l))} \\
\leq C t^{\f{m+1}{2}} \max_{\substack{ (i_1,\ldots,i_k)\in \mathcal{A} \\ m < \deg(i_1,\ldots,i_k)\leq m+2}} \sup_{r \in \mathcal{H}} |\beta^{C_1}_{i_1}\cdots\beta^{C_1}_{i_k}g(t,r)|,
\end{multline}
for $0<t<1$ and some test function $ g $ with all derivatives bounded.

\begin{remark}
Due to the a priori bounds on the moments of the solution process $ R_t $, we can estimate the probability for $ R $ to leave a
ball of radius $ C_1 $ and we can therefore control ``a priori'' the error of replacing the vector fields  $ \beta_0,\ldots,\beta_d $ by globally $ \Lip(m+1) $ vector fields $ \beta_0^{C_1},  \ldots, \beta_d^{C_1} $ on some large ball of radius $ C_1 $.
\end{remark}

For the global method (in fact an iteration due to the Markov property), divide the interval $[0,T]$ into $p$ subintervals
according to the partition $0= t_0 < t_1 < \cdots < t_p = T$. For a
multi-index $(l_1,\ldots,l_p)\in \{1,\ldots,N\}^p$ consider the path
$\omega_{l_1,\ldots,l_p}$ defined by concatenating the paths
$\omega_{l_1},\ldots,\omega_{l_p}$, i.e.~$\omega_{l_1,\ldots,l_p}(t) = \omega_
{l_1}(t)$ for $t \in ]0, t_1]$ and
\begin{equation*}
  \omega_{l_1,\ldots,l_p}(t) = \omega_{l_1,\ldots,l_p}(t_{r-1}) + 
  \omega_{l_r}(t-t_{r-1})
\end{equation*}
for $r$ such that $t\in ]t_{r-1},t_r]$, where $\omega_{l_r}$ is scaled to be a
cubature path on the interval $[0,t_r-t_{r-1}]$.
\begin{proposition}
  \label{prop:bounded_multistep}
  Fix $T>0$, $m\in \N$, $ C_1 > 0 $, a cubature formula of degree $m$ as in
  Definition~\ref{def:cubature_wiener_space} and a partition of $[0,T]$ as
  above. For every test function $ g $ there is a constant $D$ independent of the partition such
  that
  \begin{multline*}
    \sup_{r\in \mathcal{H}} \Bigl\lvert E(g(t,R^{C_1}_T)) - \sum_{(l_1,\ldots,l_p)\in
      \{1,\ldots,N\}^p} \lambda_{l_1} \cdots \lambda_{l_p}
    g(t,R^{C_1}_T(\omega_{l_1,\ldots,l_p})) \Bigr\rvert \\ \leq D T
    \max_{r=1,\ldots,p} (t_r-t_{r-1})^{(m-1)/2}.
  \end{multline*}
Additionally $ R^{C_1}_T(\omega_{l_1,\ldots,l_p}) $, due to \cite{kus04}, we can allow a local error of order $ \f{m+1}{2} $ along each $ \omega_{l_j} $.
\end{proposition}

Due to Lemma \ref{transfer-sde-spde} we can transfer the previous result including the rate of convergence on the original space. Notice that the projection of the equations with vector fields $ \beta_0^{C_1},  \ldots, \beta_d^{C_1} $ only coincide on some bounded set of the original Hilbert space $ H $ with the original equation, which is, however, for numerical purposes sufficient. The transfer works so well due to the linearity of the semigroup and the projection.

\begin{remark}
The same techniques as in \cite{BayTei:08} for the inclusion of finite activity jump processes also work in this setting. We do not outline this aspect here, since our main purpose was to show that high-order weak approximation schemes exist in the realm of SPDEs under fairly general assumptions on vector fields and test functions.
\end{remark}

We can summarize the method as follows:
\begin{itemize}
 \item Choose a degree of accuracy $ m \geq 2 $ and a set of cubature paths $ \omega_1,\ldots,\omega_N $.
 \item Choose trajectories $ \omega_{l_1,\ldots,l_p} $ by means of a MC-procedure.
 \item Calculate numerically, with error of order $ \f{m+1}{2} $, the solution of the PDE obtained by ``evaluating'' the SPDE \eqref{spde-numerics} along $ \omega_{l_j} $.
 \item Apply the main result to obtain a high order convergence scheme of order $ \f{m-1}{2} $.
\end{itemize}

\begin{remark}
The advantage of high-order schemes becomes visible when the calculation of each local step is expensive: in this case a small number $ p $ is a true advantage.
\end{remark}

\begin{appendix}

\section{Stochastic Fubini theorem with respect to Poisson measures}
\label{a}
In this appendix, we provide a stochastic Fubini theorem with
respect to compensated Poisson random measures, see Theorem
\ref{thm-Fubini-Poisson}, which we require for the proof of Lemma \ref{lemma-mild-weak}.

We could not find a proof in the literature. In the appendix of
\cite{BKR}, it is merely mentioned that it can be provided the same
way as in \cite{Protter}, where stochastic integrals with respect to
semimartingales are considered. The stochastic Fubini theorem
\cite[Thm. 5]{Applebaum-OU}, which is used in the proof of
\cite[Prop. 5.3]{Knoche2}, only deals with finite measure spaces.

We start with an auxiliary result.

\begin{lemma}\label{lemma-approx}
Let $(\Omega_i,\mathcal{F}_i,\mu_i)$, $i=1,2$ be two $\sigma$-finite
measure spaces. We define the product space
\begin{align*}
(\Omega,\mathcal{F},\mu) := (\Omega_1 \times \Omega_2, \mathcal{F}_1
\otimes \mathcal{F}_2, \mu_1 \otimes \mu_2).
\end{align*}
For each $\Phi \in L^2(\Omega, \mathcal{F}, \mu)$ there exists a
sequence
\begin{align}\label{Phi-in-span}
(\Phi_n)_{n \in \mathbb{N}} \subset {\rm span} \{ \mathbbm{1}_{A_1}
\mathbbm{1}_{A_2} : A_i \in \mathcal{F}_i \text{ with $\mu_i(A_i) <
\infty$, $i=1,2$} \}
\end{align}
such that $\Phi_n \rightarrow \Phi$ in $L^2(\Omega, \mathcal{F},
\mu)$.
\end{lemma}

\begin{proof}
Let $ \Phi \in L^2(\Omega, \mathcal{F}, \mu)$ be arbitrary. We decompose $\Phi
= \Phi^+ - \Phi^-$ into its positive and negative part. There are
sequences $(\Phi_n^+)_{n \in \mathbb{N}}$, $(\Phi_n^-)_{n \in
\mathbb{N}}$ of nonnegative  measurable functions, taking only a
finite number of values, such that $\Phi_n^+ \uparrow \Phi^+$ and
$\Phi_n^- \uparrow \Phi^-$, see, e.g., \cite[Satz 11.6]{Bauer}.

Moreover, since $\mu_1$ and $\mu_2$ are $\sigma$-finite measures,
there exist sequences $(C_n) _{n \in \mathbb{N}} \subset
\mathcal{F}_1$ and $(D_n) _{n \in \mathbb{N}} \subset \mathcal{F}_2$
such that $\mu_1(C_n) < \infty$, $\mu_2(D_n) < \infty$ for all $n
\in \mathbb{N}$ and $C_n \uparrow \Omega_1$, $D_n \uparrow \Omega_2$
as $n \rightarrow \infty$. By Lebesgue's dominated convergence
theorem we have $(\Phi_n^+ - \Phi_n^-) \mathbbm{1}_{C_n \times D_n}
\rightarrow \Phi$ in $L^2(\Omega, \mathcal{F}, \mu)$.

Therefore, we may, without loss of generality, assume that $\Phi =
\sum_{j=1}^m c_j \mathbbm{1}_{A_j}$, where $m \in \mathbb{N}$, $c_j
\in \mathbb{R} \setminus \{ 0 \}$, $j = 1,\ldots,m$ and $A_j \in
(\mathcal{F}_1 \otimes \mathcal{F}_2) \cap (C_1 \times C_2)$, $j =
1,\ldots,m$, where $C_i \in \mathcal{F}_i$, $i=1,2$ and $\mu_i(C_i)
< \infty$, $i = 1,2$.

Note that the trace $\sigma$-algebra $(\mathcal{F}_1 \otimes
\mathcal{F}_2) \cap (C_1 \times C_2)$ is generated by the algebra
\begin{align*}
\mathcal{A} = \bigg\{ \biguplus_{k=1}^p D_k \times E_k \,|\, p \in
\mathbb{N} \text{ and $D_k \in \mathcal{F}_1 \cap C_1$, $E_k \in
\mathcal{F}_2 \cap C_2$ for $k = 1,\ldots,p$} \bigg\}.
\end{align*}
By \cite[Satz 5.7]{Bauer} there exists, for each $j \in \{
1,\ldots,m \}$ and each $n \in \mathbb{N}$, a set $B_j^n \in
\mathcal{A}$ such that $\mu(A_j \Delta B_j^n) < \frac{1}{m^2 n
c_j^2}$.

Setting $\Phi_n := \sum_{j=1}^m c_j \mathbbm{1}_{B_j^n}$ for $n \in
\mathbb{N}$ we have (\ref{Phi-in-span}) and
\begin{align*}
&\int_{\Omega} |\Phi(\omega) - \Phi_n(\omega)|^2 d \mu(\omega) \leq
\int_{\Omega} \bigg( \sum_{j=1}^m |c_j| \cdot
|\mathbbm{1}_{A_j}(\omega) - \mathbbm{1}_{B_j^n}(\omega)| \bigg)^2 d
\mu(\omega)
\\ &\leq m \sum_{j=1}^m \int_{\Omega} c_j^2
\mathbbm{1}_{A_j \Delta B_j^n}(\omega) d \mu(\omega) = m
\sum_{j=1}^m c_j^2 \mu(A_j \Delta B_j^n) < \frac{1}{n}, \quad n \in
\mathbb{N}
\end{align*}
showing that $\Phi_n \rightarrow \Phi$ in $L^2(\Omega, \mathcal{F},
\mu)$.
\end{proof}

Let $T \in \mathbb{R}_+$ be a finite time horizon. In order to have
a more convenient notation in the following stochastic Fubini
theorem, we introduce the spaces
\begin{align*}
L_T^2(\mu) &:= L_T^2(\mu;\mathbb{R}),
\\ L_T^2(\lambda) &:= L_T^2([0,T],\mathcal{B}[0,T],\lambda),
\\ L_T^2(\mathbb{P} \otimes \lambda) &:= L^2(\Omega \times [0,T],\mathcal{F}_T \otimes
\mathcal{B}[0,T], \mathbb{P} \otimes \lambda),
\end{align*}
where $L_T^2(\mu;H)$ for a separable Hilbert space $H$ was defined
in \eqref{integrands-poisson_square_integrable}, and
\begin{align}\label{def-prod-mu}
L_T^p(\mu \otimes \lambda) &:= L^p(\Omega \times
[0,T] \times E \times [0,T], \mathcal{P}_T \otimes \mathcal{E}
\otimes \mathcal{B}[0,T], \mathbb{P} \otimes \lambda \otimes F
\otimes \lambda)
\end{align}
for all $p \geq 1$.

\begin{theorem}\label{thm-Fubini-Poisson}
For each $\Phi \in L_T^2(\mu \otimes \lambda)$ we have
\begin{align}\label{int1-in-Fubini}
\int_0^T \Phi(\cdot,\cdot,s) ds \in L_T^2(\mu),
\end{align}
there exists $\phi \in L_T^2(\mathbb{P} \otimes \lambda)$ such
that for $\lambda$-almost all $s \in [0,T]$
\begin{align}\label{int2-in-Fubini}
\phi(s) = \int_0^T \int_E \Phi(t,x,s) (\mu(dt,dx) - F(dx) dt) \quad
\text{in $L^2(\Omega,\mathcal{F}_T,\mathbb{P})$}
\end{align}
and we have the identity
\begin{align}\label{Fubini-formula}
\int_0^T \phi(s) ds = \int_0^T \int_E \bigg( \int_0^T \Phi(t,x,s) ds
\bigg) (\mu(dt,dx) - F(dx) dt) \quad \text{in
$L^2(\Omega,\mathcal{F}_T,\mathbb{P})$.}
\end{align}
\end{theorem}

\begin{proof}
Let $V \subset L_T^2(\mu \otimes \lambda)$ be the vector space
\begin{align*}
V := {\rm span} \{ Kf \, | \, K \in L_T^2(\mu), f \in L_T^2(\lambda)
\}.
\end{align*}
Let $\Phi \in V$ be arbitrary. Then there exist $n \in \mathbb{N}$
and $c_i \in \mathbb{R}$, $K_i \in L_T^2(\mu)$, $f_i \in
L_T^2(\lambda)$, $i = 1,\ldots,n$ such that $\Phi = \sum_{i=1}^n c_i
K_i f_i$. Moreover we have
\begin{align*}
\phi &:= \int_0^T \int_E \Phi(t,x,\cdot) (\mu(dt,dx) - F(dx) dt)
\\ &= \sum_{i=1}^n c_i f_i(\cdot) \int_0^T \int_E K_i(t,x) (\mu(dt,dx) -
F(dx) dt) \in L_T^2(\mathbb{P} \otimes \lambda),
\\ &\int_0^T \Phi(\cdot,\cdot,s) ds = \sum_{i=1}^n c_i K_i(\cdot,\cdot) \int_0^T
f_i(s) ds \in L_T^2(\mu)
\end{align*}
and identity (\ref{Fubini-formula}) is valid.

For each $\Phi \in L_T^2(\mu \otimes \lambda) \cap L_T^1(\mu \otimes
\lambda)$ we have, according to \cite[Prop. II.1.14]{JS},
\begin{equation}\label{split-integral}
\begin{aligned}
&\int_0^T \int_E \Phi(t,x,\cdot) (\mu(dt,dx) - F(dx) dt)
\\ &= \sum_{n \in \mathbb{N}} \Phi(\tau_n,\beta_{\tau_n},\cdot)
\mathbbm{1}_{\{ \tau_n \leq T \}} - \int_0^T \int_E \Phi(t,x,\cdot)
F(dx) dt,
\end{aligned}
\end{equation}
where $(\tau_n)_{n \in \mathbb{N}}$ is a sequence of stopping times
and $\beta$ denotes an $E$-valued optional process. By the classical
Fubini theorem we deduce that the stochastic integral in
(\ref{split-integral}) is $\mathcal{F}_T \otimes
\mathcal{B}[0,T]$-measurable. Using the It\^o-isometry
(\ref{Ito-isometry-mu}) we obtain
\begin{align*}
&\int_0^T \mathbb{E} \left[ \bigg( \int_0^T \int_E \Phi(t,x,s)
(\mu(dt,dx) - F(dx) dt) \bigg)^2 \right] ds
\\ &= \int_0^T \mathbb{E} \bigg[ \int_0^T \int_E |\Phi(t,x,s)|^2 F(dx) dt \bigg] ds
< \infty,
\end{align*}
because $\Phi \in L_T^2(\mu \otimes \lambda)$ by hypothesis, and we
conclude
\begin{align}\label{integral-measurable}
\int_0^T \int_E \Phi(t,x,\cdot) (\mu(dt,dx) - F(dx) dt) \in
L_T^2(\mathbb{P} \otimes \lambda), \quad \Phi \in L_T^2(\mu \otimes
\lambda) \cap L_T^1(\mu \otimes \lambda).
\end{align}
Now let $\Phi \in L_T^2(\mu \otimes \lambda)$ be arbitrary. By the
classical Fubini theorem the integral appearing in
(\ref{int1-in-Fubini}) is $\mathcal{P}_T \otimes
\mathcal{E}$-measurable. H\"older's inequality and the hypothesis
$\Phi \in L_T^2(\mu \otimes \lambda)$ yield
\begin{align*}
&\mathbb{E}\left[ \int_0^T \int_E \bigg( \int_0^T \Phi(t,x,s) ds
\bigg)^2 F(dx) dt \right]
\\ &\leq T \mathbb{E}\bigg[ \int_0^T
\int_E \int_0^T |\Phi(t,x,s)|^2 ds F(dx) dt \bigg] < \infty,
\end{align*}
and hence (\ref{int1-in-Fubini}) is valid.

Since the measure $F$ is $\sigma$-finite, there exists a sequence
$(B_n)_{n \in \mathbb{N}} \subset E$ with $F(B_n) < \infty$, $n \in
\mathbb{N}$ and $B_n \uparrow E$. We define
\begin{align*}
\phi_n := \int_0^T \int_E \Phi(t,x,\cdot) \mathbbm{1}_{B_n}(x)
(\mu(dt,dx) - F(dx) dt), \quad n \in \mathbb{N}.
\end{align*}
By (\ref{integral-measurable}) we have $\phi_n \in L_T^2(\mathbb{P}
\otimes \lambda)$ for all $n \in \mathbb{N}$. Now, we shall prove
that $(\phi_n)_{n \in \mathbb{N}}$ is a Cauchy sequence in
$L_T^2(\mathbb{P} \otimes \lambda)$.

Let $\epsilon > 0$ be arbitrary. By Lebesgue's theorem, there exists
an index $n_0 \in \mathbb{N}$ such that
\begin{align*}
\int_0^T \mathbb{E} \bigg[ \int_0^T \int_E | \Phi(t,x,s) |^2
\mathbbm{1}_{E \setminus B_n}(x) F(dx) dt \bigg] ds < \epsilon,
\quad n \geq n_0
\end{align*}
For all $m > n \geq n_0$ we obtain by the It\^o-isometry
(\ref{Ito-isometry-mu})
\begin{align*}
&\int_0^T \mathbb{E} [| \phi_n(s) - \phi_m(s) |^2] ds = \int_0^T
\mathbb{E} \bigg[ \int_0^T \int_E | \Phi(t,x,s) |^2 \mathbbm{1}_{B_m
\setminus B_n}(x) F(dx) dt \bigg] ds
\\ &\leq \int_0^T
\mathbb{E} \bigg[ \int_0^T \int_E | \Phi(t,x,s) |^2 \mathbbm{1}_{E
\setminus B_n}(x) F(dx) dt \bigg] ds < \epsilon,
\end{align*}
establishing that $(\phi_n)_{n \in \mathbb{N}}$ is a Cauchy sequence
in $L_T^2(\mathbb{P} \otimes \lambda)$. Thus, there exists $\phi \in
L_T^2(\mathbb{P} \otimes \lambda)$ such that $\phi_n \rightarrow
\phi$ in $L_T^2(\mathbb{P} \otimes \lambda)$. The relation
\begin{align*}
\int_0^T \mathbb{E}[| \phi_n(s) - \phi(s) |^2]ds \rightarrow 0 \quad
\text{for $n \rightarrow \infty$}
\end{align*}
implies that there exists a subsequence $(n_k)_{k \in \mathbb{N}}$
such that
\begin{align*}
\mathbb{E}[| \phi_{n_k}(s) - \phi(s) |^2] \rightarrow 0 \quad
\text{for $\lambda$-almost all $s \in [0,T]$,}
\end{align*}
that is $\phi_{n_k}(s) \rightarrow \phi(s)$ in
$L_T^2(\Omega,\mathcal{F}_T,\mathbb{P})$ for $\lambda$-almost all $s
\in [0,T]$. We define
\begin{align*}
\psi := \int_0^T \int_E \Phi(t,x,\cdot) (\mu(dt,dx) - F(dx) dt).
\end{align*}
By the classical Fubini theorem we have $\Phi(\cdot,\cdot,s) \in
L_T^2(\mu)$ for $\lambda$-almost all $s \in [0,T]$. The
It\^o-isometry (\ref{Ito-isometry-mu}) and Lebesgue's theorem yield
\begin{align*}
\mathbb{E}[|\psi(s) - \phi_n(s)|^2] = \mathbb{E}\bigg[ \int_0^T
\int_E |\Phi(t,x,s)|^2 \mathbbm{1}_{E \setminus B_n}(x) F(dx) dt
\bigg] \rightarrow 0 \quad \text{for $n \rightarrow \infty$,}
\end{align*}
implying $\phi_n(s) \rightarrow \psi(s)$ in
$L_T^2(\Omega,\mathcal{F}_T,\mathbb{P})$ for $\lambda$-almost all $s
\in [0,T]$. We infer that $\phi(s) = \psi(s)$ in
$L_T^2(\Omega,\mathcal{F}_T,\mathbb{P})$ for $\lambda$-almost all $s
\in [0,T]$, proving (\ref{int2-in-Fubini}).

According to Lemma \ref{lemma-approx} there exists a sequence
$(\Phi_n)_{n \in \mathbb{N}} \subset V$ such that $\Phi_n
\rightarrow \Phi$ in $L_T^2(\mu \otimes \lambda)$. From the
beginning of the proof we know that for each $n \in \mathbb{N}$ we
have
\begin{align*}
\int_0^T \Phi_n(\cdot,\cdot,s) ds \in L_T^2(\mu), \quad \int_0^T
\int_E \Phi_n(t,x,\cdot) (\mu(dt,dx) - F(dx) dt) \in
L_T^2(\mathbb{P} \otimes \lambda)
\end{align*}
and the identity
\begin{equation}\label{Fubini-formula-n}
\begin{aligned}
&\int_0^T \bigg( \int_0^T \int_E \Phi_n(t,x,s) (\mu(dt,dx) - F(dx)
dt) \bigg) ds
\\ &= \int_0^T \int_E \bigg( \int_0^T \Phi_n(t,x,s) ds
\bigg) (\mu(dt,dx) - F(dx) dt)
\end{aligned}
\end{equation}
in $L_T^2(\Omega,\mathcal{F}_T,\mathbb{P})$. By H\"older's
inequality, (\ref{int2-in-Fubini}), the It\^o-isometry
(\ref{Ito-isometry-mu}) and the convergence $\Phi_n \rightarrow
\Phi$ in $L_T^2(\mu \otimes \lambda)$ we get
\begin{equation}\label{Fubini-conv-1}
\begin{aligned}
&\mathbb{E}\left[ \bigg( \int_0^T \bigg( \int_0^T \int_E
\Phi_n(t,x,s) (\mu(dt,dx) - F(dx) dt) \bigg) ds - \int_0^T \phi(s)
ds \bigg)^2 \right]
\\ &\leq T \int_0^T \mathbb{E} \left[ \bigg( \int_0^T \int_E \Phi_n(t,x,s)
(\mu(dt,dx) - F(dx) dt) - \phi(s) \bigg)^2 \right] ds
\\ &= T \int_0^T \mathbb{E} \bigg[ \int_0^T \int_E | \Phi_n(t,x,s)
- \Phi(t,x,s) |^2 F(dx) dt \bigg] ds \rightarrow 0.
\end{aligned}
\end{equation}
The It\^o-isometry (\ref{Ito-isometry-mu}), H\"older's inequality
and the convergence $\Phi_n \rightarrow \Phi$ in $L_T^2(\mu \otimes
\lambda)$ yield
\begin{equation}\label{Fubini-conv-2}
\begin{aligned}
&\mathbb{E} \Bigg[ \bigg( \int_0^T \int_E \bigg( \int_0^T
\Phi_n(t,x,s) ds \bigg) (\mu(dt,dx) - F(dx) dt)
\\ &\quad - \int_0^T
\int_E \bigg( \int_0^T \Phi(t,x,s) ds \bigg) (\mu(dt,dx) - F(dx) dt)
\bigg)^2 \Bigg]
\\ &= \mathbb{E} \left[ \int_0^T \int_E \bigg( \int_0^T (\Phi_n(t,x,s) - \Phi(t,x,s)) ds \bigg)^2 F(dx) dt \right]
\\ &\leq T \mathbb{E} \bigg[ \int_0^T \int_E \int_0^T |\Phi_n(t,x,s) - \Phi(t,x,s)|^2 F(dx) dt ds
\bigg] \rightarrow 0.
\end{aligned}
\end{equation}
Combining (\ref{Fubini-formula-n}), (\ref{Fubini-conv-1}) and
(\ref{Fubini-conv-2}) we arrive at (\ref{Fubini-formula}).
\end{proof}

\end{appendix}

\end{document}